\documentclass{article}
\usepackage{etex}
\usepackage[latin1]{inputenc}
\usepackage{amsmath, enumerate, amsfonts}
\usepackage[widemargins]{a4}
\usepackage{tikz} 
\usepackage{tikz-cd}  
\usepackage{amsmath, enumerate, amsfonts,amssymb, pictexwd, amscd}

\newcommand{\pcom}{^\wedge_p}
\usepackage[all]{xy}
\xyoption{arc}
\newtheorem{Conjecture}{Conjecture}[section]

\newtheorem{Theorem}{Theorem}[section]

\newtheorem{Proposition}{Proposition}[section]
\newtheorem{Lemma}{Lemma}[section]

\newcommand*{\F}{\mathbb{F}}

\newcommand*{\hoTop}{\textbf{hoTop}}
\DeclareMathOperator{\Aut}{Aut}
\DeclareMathOperator{\ord}{ord}
\DeclareMathOperator{\Ker}{Ker}
\DeclareMathOperator{\Inn}{Inn}
\DeclareMathOperator{\id}{id}
\DeclareMathOperator{\Out}{Out}
\DeclareMathOperator{\Hom}{Hom}
\DeclareMathOperator{\GL}{GL}
\DeclareMathOperator{\PSL}{PSL}
\DeclareMathOperator{\Mor}{Mor}
\DeclareMathOperator{\hocolim}{hocolim}
\title{Signalizer functors, existence, and the fundamental group}
\author{Nora Seeliger}
\begin{document}
\maketitle
\begin{abstract}
 We solve the seventh problem of Oliver's list [M.\ Aschbacher, R.\ Kessar, B.\ Oliver, \textit{Fusion systems in algebra and topology}, LMS Lecture Note Series: 31, Cambridge University Press, 2011] via an explicit signalizer functor construction in the sense of Aschbacher-Chermak for various group models. Moreover we prove the existence of centric linking systems via group models in certain cases which is the first problem  and give applications to the fundamental group which is the eighth problem of the list respectively. We illustrate with many examples.
 \end{abstract}
\tableofcontents
\section{Introduction}
\label{cIntro}
Motivated by the Martino-Priddy Conjecture stated in \cite{MartinoPriddy} and proved by Oliver in \cite{MP1}, \cite{MP2} the theory of $p-$local finite groups introduced by Broto,
 Levi and Oliver \cite{BLO2} describes spaces which are similar to $p-$completed classifying spaces of finite groups via 
  $p-$local structure, at least up to $\mathbb{F}_p$-cohomology. 
   A $p$-local finite group is a
   triple $(S,\mathcal{F},\mathcal{L})$
    which encodes many properties similar to those encoded in a finite group and every finite group gives rise to a $p$-local finite group.
However this is not true for all of them and a $p$-local finite group is not a group. It is therefore an interesting and natural question to ask whether it is possible to find any group which encodes all the structure of a given $p$-local finite group and such that the $p$-local finite group can be recovered from it. This means a
group which has the appropriate Sylow $p$-subgroup and is model for the fusion system and also a group model for the
linking system induced by an explicit signalizer functor, and in addition that the classifying space has $\mathbb{F}_p$-cohomology isomorphic to the stable elements.  Every $p$-local finite group which is induced by a finite group has such a group model provided by this finite group as shown by Broto, Levi, Oliver \cite{BLO2}. The general case remains open.  The fact that the linking system can be simply connected shows that it cannot always be equivalent to the classifying space of a group model however there can still be a group model whose classifying space has cohomology isomorphic to the stable elements.
 The world of groups
 is big which is why we hope that such a group model could exists even though no explicit construction proves this so far.\\
  We start by giving a review of various previous results putting our own work into perspective. Kan and Thurston state \cite{KanThurston} that for every connected space there is a map from a $K(\pi,1)$ which is an isomorphism in homology for any system of local coefficients. The group they construct however does not realize neither the Sylow $p$-subgroup, nor the fusion system nor the linking system nor has a signalizer functor.
   The first infinite group model for a saturated fusion system was constructed by Aschbacher and Chermak \cite{AschbacherChermak}
realizing the Solomon groups \cite{LO}. They also gave the only
explicit description of a signalizer functor for a group model
    so far. The first conceptual treatment of group models for fusion systems are due to Leary-Stancu \cite{Ian+Radu} and Robinson \cite{Robinson1}. These groups realize the linking system as well induced by a signalizer functor as proved in \cite{LS} (where the signalizer functor is not constructed explicitely) however
     both constructions
     fail to have the right cohomology in the general case \cite[Proposition
     4.2]
     {gmffs} and Proposition
      \ref{learystancupbad}
   this article.\\ 
   In this note  we construct a new group model for saturated fusion systems involving amalgams
   to generalize and overcome the shortcomings of previous constructions.
   (In Proposition \ref{learystancupbad} we discuss why the method of iterated HNN-constructions seems to be less promising for our purposes.)
   \\ Moreover
    we construct an explicit signalizer functor for the Robinson and our new
     group model.
     This is a more direct and most
    important complete  solution
      to the seventh problem of Oliver's list \cite{AKO}
    in the sense that it provides an explicit formula
     than the second author's previous work with Libman \cite{LS}. We
      believe that independently from our context signalizer functor constructions are of great interest to group theorists.
The continued study of group models for fusion systems is motivated by the fact that a group model which has the cohomology of the saturated fusion system provides an independent proof of the existence of a centric linking system for this saturated fusion system. This is a much stronger version of the existence of centric linking systems
which is a theorem proved by Chermak \cite{Chermak} and Oliver \cite{Oliver}
 using the classification of finite simple groups. In addition
 we give some applications to the fundamental group. The question whether for every $p$-local finite group there exists a group model which has the $\mathbb{F}_p$-cohomology isomorphic to the stable elements remains and will be addressed in the last subsection. We prove this conjecture in some special cases.\\
The author met Malte Leip at the Max-Planck-Institut
f\"ur Mathematik in Bonn in August 2012. He contributed results in some special cases and thanks the MPI for their hospitality during that month, and the author for introducing him to this topic.
The author was supported by ANR BLAN08-2-338236, an
Erwin-Schr\"odinger-Institute Junior-Research-Fellowship, a
Mathematisches-Forschungsinstitut-Oberwolfach Leibniz-Fellowship, an
invitation to the Max-Planck-Institut f\"ur Mathematik in Bonn, the Centre Recerca Matem\`atica Barcelona, the Australian Research Council Discovery Project DP120101399 at the ANU, and a fellowship from the University of Haifa.

\section{Preliminaries}
\subsection{Fusion Systems}
We review basic definitions of fusion systems and
linking systems and establish our notations. Our main references are
\cite{BLO2}, \cite{BLO4} and \cite{IntroMarkus}. Let $S$ be a finite
$p$-group. A \textbf{fusion system} $\mathcal{F}$ on $S$ is a
category whose objects are all the subgroups of $S$, and which
satisfies the following two properties for all $P,Q\leq S$: The set
$Hom_{\mathcal{F}}(P,Q)$ contains injective group homomorphisms and
amongst them all morphisms induced by conjugation of elements in $S$
and each element is the composite of an isomorphism in $\mathcal{F}$
followed by an inclusion. Two subgroups $P,Q\leq S$ will be called
\textbf{$\mathcal{F}-$conjugate} if they are isomorphic in
$\mathcal{F}$. Define
$Out_{\mathcal{F}}(P)=Aut_{\mathcal{F}}(P)/Inn(P)$ for all $P\leq
S$.
 A subgroup $P\leq S$ is \textbf{fully centralized} resp.\ \textbf{fully normalized} in $\mathcal{F}$ if $|C_S(P)|\geq |C_S(P')|$ resp.\ $|N_S(P)|\geq |N_S(P')|$ for all $P'\leq S$ which is $\mathcal{F}$-conjugate to $P$.
 $\mathcal{F}$ is called \textbf{saturated} if for all $P\leq S$ which is fully normalized in $\mathcal{F}$, $P$ is fully centralized in $\mathcal{F}$ and $Aut_S(P)\in Syl_p(Aut_{\mathcal{F}}(P))$ and moreover if $P\leq S$ and $\phi\in Hom_{\mathcal{F}}(P,S)$ are such that $\phi (P)$ is fully centralized, and if we set
$N_{\phi}=\{g\in N_S(P)|\phi c_g \phi^{-1}\in Aut_S(\phi (P))\}$,
then there is $\overline{\phi}\in Hom_{\mathcal{F}}(N_{\phi},S)$
such that $\overline{\phi}|_P=\phi$. A subgroup $P\leq S$ will be
called $\mathcal{F}$-\textbf{centric} if $C_S(P')\leq P'$ for all
$P'$ which are $\mathcal{F}-$conjugate to $P$. Denote
$\mathcal{F}^c$ the full subcategory of $\mathcal{F}$ with objects
the $\mathcal{F}$-centric subgroups of $S$. 
Let
$\mathcal{O}(\mathcal{F})$ be the \textbf{orbit category} of
$\mathcal{F}$ with objects the same objects as $\mathcal{F}$ and morphisms the set
 $Mor_{\mathcal{O}(\mathcal{F})}(P,Q)=Mor_{\mathcal{F}}(P,Q)/Inn(Q)$.
 Let $\mathcal{O}^c(\mathcal{F})$ be the full subcategory of
 $\mathcal{O}(\mathcal{F})$ with objects the $\mathcal{F}-$centric
 subgroups of $\mathcal{F}$. 
 Let $G$ be a discrete group. A finite subgroup $S$ of $G$ will be
called a \textbf{Sylow} $p$-subgroup of $G$ if $S$ is a $p$-subgroup
of $G$ and all $p-$subgroups of $G$ are conjugate to a subgroup of $S$. A group $G$ is called $p$\textbf{-perfect} if it has no normal subgroup of index $p$. Let $\mathcal{F}$ be a saturated fusion
system over the finite $p$-group $S$. A subgroup $P\leq S$ is called
$\mathcal{F}$-\textbf{radical} if
$O_p(Aut_{\mathcal{F}}(P))=Aut_P(P)$, where $O_p(-)$ denotes the
maximal normal $p-$subgroup. A subgroup $P\leq S$ which is both $\mathcal{F}$-centric and $\mathcal{F}$-radical is called $\mathcal{F}$-centric-radical. Denote
$\mathcal{F}^{cr}$ the full subcategory of $\mathcal{F}$ with objects
the $\mathcal{F}$-centric-radical subgroups of $S$. A subgroup $P\leq S$ is called
$\mathcal{F}$\textbf{-essential} if P is $\mathcal{F}$-centric and
$Aut_{\mathcal{F}}(P)/Aut_P(P)$ has a strongly $p-$embedded
subgroup, i. e. there exists a proper subgroup $Q<P$ which contains
Sylow $p-$subgroup $R$ of $P$ and $R\neq 1$ but $Q\cap ^xR=1$ for any
$x\in P-Q$. A subgroup $Q\triangleleft S$ is \textbf{normal} in $\mathcal{F}$ if each $\alpha \in Hom_{\mathcal{F}(P,P')}$ extends to a morphism $\overline{\alpha}\in Hom_{\mathcal{F}}(PQ,P'Q)$ which sends $Q$ to itself. A saturated fusion system is called \textbf{constrained} if there is some $Q\triangleleft S$ which is $\mathcal{F}$-centric and normal in $\mathcal{F}$.\\
Let $\mathcal{F}$ be a fusion system over a finite $p$-group $S$. A saturated fusion subsystem of $\mathcal{F}$ is a subcategory $\mathcal{E}\subseteq \mathcal{F}$ which is itself a saturated fusion system over a subgroup $T\leq S$.
A fusion subsystem $\mathcal{E}\subseteq\mathcal{F}$ over $T\unlhd S$ is $\mathcal{F}$-invariant if $T$ is strongly closed in $\mathcal{F}$, and 
 $^{\alpha}\mathcal{E} =\mathcal{E}$ for each $\alpha\in Aut_{\mathcal{F}}(T)$, and
 for each $P\leq T$ and each $\phi \in Hom_{\mathcal{F}}
(P,T)$,
there are 
$\alpha \in Aut_{\mathcal{F}}
(T)$
 and $\phi _0\in Hom_{\mathcal{E}}(P,T)$ such that $\phi =\alpha \circ \phi _0$.
A
fusion subsystem $\mathcal{E}\subseteq\mathcal{F}$ over $T\unlhd S$ is
\textbf{normal} 
$\mathcal{E}\unlhd\mathcal{F}$ if $\mathcal{E}$ is \textbf{weakly normal}, and each $\alpha\in Aut_{\mathcal{E}}(T)$ 
extends to some $\overline{\alpha}\in Aut_{\mathcal{F}}(TC_S(T))$ such that
$[\overline{\alpha},C_S(T)]\leq Z(T)$. The fusion system is \textbf{simple} if it contains no proper nontrivial normal fusion subsystem.
Let $S$ be a finite $p-$group and let $P_1,...,P_r,Q_1,...,Q_r$ be
subgroups of $S$. Let $\phi _1,...,\phi _r$ be injective group
homomorphisms $\phi _i:P_i\rightarrow Q_i$ $\forall i$. The fusion
system \textbf{generated} by $\phi _1,...,\phi_r$ is the minimal
fusion system $\mathcal{F}$ over $S$ containing $\phi _1,...,\phi
_r$. Let $P_1,\dots,P_n$ be a collection of $\mathcal{F}$-centric
subgroups of $S$, which might contain multiples.
Let $K_1,\dots K_n$ be groups with $K_i$ is a subgroup of $\Aut_{\mathcal{L}}(P_i)$
containing $\delta (P_i)$ for all $i=1,\cdots ,n$. 
We say that $K=\{K_1,\dots K_n\}$ is  \textbf{generating} if $\pi(K_1),\dots,\pi(K_n)$ generate
$\mathcal{F}$. Let $\mathcal{F}$ be a fusion system on a finite
$p-$group $S$. A subgroup $T\leq S$ is \textbf{strongly closed} in
$S$ with respect to $\mathcal{F}$, if for each subgroup $P$ of $T$,
each $Q\leq S$, and each $\phi\in Mor_{\mathcal{F}}(P,Q)$, $\phi
(P)\leq T$.
Fix any pair $S\leq G$, where $G$ is a (possibly infinite) group and $S$ is a finite $p-$subgroup.
    Define $\mathcal{F}_S(G)$ to be the category whose objects are the subgroups of $S$, and where
    $Mor_{\mathcal{F}_S(G)}(P,Q)=Hom_G(P,Q)=\{c_g\in Hom(P,Q)|g\in G, gPg^{-1}\leq Q\}
    \cong N_G(P,Q)/C_G(P).$
    Here $c_g$ denotes the homomorphism \textbf{conjugation} by $g$ $(x\mapsto gxg^{-1})$.
For each $P\leq S$, let $C'_G(P)$ be the maximal $p-$perfect subgroup of $C_G(P)$.Let $\mathcal{F}$ be a fusion system over a $p$-group $S$. For a discrete group $G$ and any finite set $\mathcal{H}$ of finite subgroups of $G$, let $\mathcal{T}_{\mathcal{H}}(G)$ denote the \textbf{transporter category of $G$}: the category with $Ob(\mathcal{T}_{\mathcal{H}}(G))=\mathcal{H}$, and where for each $P,Q\in\mathcal{H}$, $Mor_{\mathcal{T}_{\mathcal{H}}(G)}(P,Q)=N_G(P,Q)=\{g\in G|gPg^{-1}\leq Q\}$ the \textbf{transporter set}. The category  $\mathcal{T}^c_{S}(G)$ is the full subcategory of $\mathcal{T}_{S}(G)$
  with objects the $\mathcal{F}_S(G)$-centric subgroups of $S$.
 A \textbf{linking system} $\mathcal{L}$
associated with $\mathcal{F}$ is a finite category together
with a pair of functors $\mathcal{T}_{Ob(\mathcal{L})}(S)\overset{\delta}{\longrightarrow}\mathcal{L}\overset{\pi }{\longrightarrow}\mathcal{F}$ such that the following conditions are satisfied: $Ob(\mathcal{L})$ is a set of subgroups of $S$ closed under $\mathcal{F}$-conjugacy and overgroups, and includes all subgroups which are $\mathcal{F}$-centric and $\mathcal{F}$-radical. Also $\delta$ is the
identity on objects, and $\pi$ is the inclusion on objects. and surjective on morphisms. For
each pair of objects $P,Q\in ob(\mathcal{L})$ such that $P$ is fully centralized in $\mathcal{F}$, $C_S(P)$ acts freely on
$Mor_{\mathcal{L}}(P,Q)$ via $\delta _{P,P}$ by right composition and $\pi _{P,Q}$ induces
a bijection
$Mor_{\mathcal{L}}(P,Q)/C_S(P)\overset{\simeq}{\longrightarrow}Hom_{\mathcal{F}}(P,Q).$
For each $P,Q\in ob(\mathcal{L})$ and each $g\in N_S(P,Q)$,
$\pi _{P,Q}$ sends $\delta _{P,Q}\in Mor_{\mathcal{L}}(P,Q)$ to $=c_g \in Hom_{\mathcal{F}}(P,Q)$. For all $\phi\in
Mor_{\mathcal{L}}(P,Q)$ and all $g\in P, \phi\circ\delta_{P,P}(g)=\delta
_{Q,Q}(\pi (\phi (g)))\circ \phi$.
A \textbf{centric linking system}
associated with a saturated fusion system $\mathcal{F}$ is a linking system with objects the set of $\mathcal{F}$-centric subgroups of $S$.
For every saturated fusion system there exists one associated centric linking system up to isomorphism \cite{Chermak}. 
 Let $\mathcal{L}^c_S(G)$ be the category whose objects are the $\mathcal{F}_S(G)-$centric subgroups of $S$, and where $Mor_{\mathcal{L}^c_S(G)}(P,Q)=N_G(P,Q)/C_G'(P).$
    Let $\pi:\mathcal{L}^c_S(G)\rightarrow\mathcal{F}_S(G)$ be the functor which is the inclusion on objects and sends the class of $g\in N_G(P,Q)$ to conjugation by $g$. For each $\mathcal{F}_S(G)-$centric subgroup $P\leq G$, let $\delta _P:P\rightarrow Aut_{\mathcal{L}^c_S(G)}(P)$ be the monomorphism induced by the inclusion $P\leq N_G(P)$.
A triple $(S,\mathcal{F},\mathcal{L})$ where $S$ is a finite
$p-$group, $\mathcal{F}$ is a saturated fusion system on $S$, and
$\mathcal{L}$ is an associated centric linking system with
$\mathcal{F}$, is called a $p-$local finite group. Its
\textbf{classifying space} is $|\mathcal{L}|\pcom$ where $(-)\pcom$
denotes the $p-$completion functor in the sense of Bousfield and
Kan. A space X is called \textbf{$p-$good} if the natural map $H_*(X;\mathbb{F}_p)\rightarrow H_*(X\pcom ;\mathbb{F}_p)$ is an isomorphism. Examples of spaces which are $p-$good are classifying spaces of finite groups. A finite group $G$ gives canonically rise to a $p-$local
finite group $(S,\mathcal{F}_S(G),\mathcal{L}^c_S(G))$ and
$BG\pcom\simeq |\mathcal{L}^c_S(G)|\pcom$ \cite{BK}. In particular, all
fusion systems coming from finite groups are saturated. In analogy with the Cartan-Eilenberg theorem the cohomology of a saturated fusion system is defined to be the inverse limit over the orbit category. Recall \cite{BLO2} there exists an isomorphism of unstable algebras between the cohomology of the classifying space of a $p-$local finite group and the cohomology of the fusion system.
Let $\mathcal{F}$ be a fusion system on the the finite $p-$group
$S$. $\mathcal{F}$ is called \textbf{Alperin} fusion system if
there are subgroups $P_1,P_2,\cdots P_r$ of $S$ and finite groups
$L_1,\cdots ,L_r$ such that for each $i$, $N_S(P_i)\in Syl_p(L_i)$,
$\mathcal{F}_{N_S(P_i)}(L_i)$ is contained in $\mathcal{F}$ and
$\mathcal{F}$ is generated by all the $\mathcal{F}_{N_S(P_i)}(L_i)$. Every saturated fusion system is Alperin
\cite[Section 4]{BLO1}. The groups $L_i=Aut_{\mathcal{L}}(P
_i)$ 
are known and unique up to isomorphism for $i=1,\cdots ,n$ regardless of the existence of $\mathcal{L}$ and also denoted $L_{P_i}$.
\\
A ring homomorphism $f:A\rightarrow B$ is called
$F-$\textbf{monomorphism} \cite{Quillen} if
every element in the kernel is nilpotent  and
$F-$\textbf{epimorphism} if every element in the cokernel is
nilpotent and $F-$\textbf{isomorphism} if it is
$F-$\textbf{monomorphism} and $F-$\textbf{epimorphism}.\\
\subsection{Group models for fusion systems}
 A saturated fusion system $\mathcal{F}$ on a finite $p-$group $S$ for which there exists no finite group $G$ such that $\mathcal{F}_S=\mathcal{F}_S(G)$ is called \textbf{exotic} (see \cite{BLO2}, chapter 9 for example). 
 A discrete group $\mathcal{G}$ is a \textbf{group model} for
$\mathcal{F}$ if $S$ is a Sylow $p-$subgroup of $\mathcal{G}$ and
$\mathcal{F}_S(\mathcal{G})=\mathcal{F}$ in the sense that all morphisms of $\mathcal{F}$ are induced by conjugation of elements in $\mathcal{G}$. Contrary to the finite case infinite groups need not to have Sylow $p$-subgroups. For constrained fusion systems there exists finite group models \cite{BCGLO1}. 
We review all general constructions for group models for fusion systems known so far \cite{Ian+Radu}, \cite{Robinson1}, \cite{LS}.

The group model constructed by Leary and Stancu is an iterated HNN-construction.
\begin{Theorem}[\cite{Ian+Radu}, Theorem 2]
Let $\mathcal{F}$ be a fusion system on $S$ generated by $\Phi=\{\phi_1, \cdots, \phi_r\}$ with $\phi _i:P_i\rightarrow Q_i$ a morphism in $\mathcal{F}$ for $P_i,Q_i$ subgroups of $S$ for $i=1,\cdots , r$. Let $T$ be a free group with free generators $t_1, \ldots, t_r$, and define $\pi_{LS}$ as the quotient of the free product $S*T$ by the relations $t_i^{-1}ut_i=\phi_i(u)$ for all $i$ and for all $u\in P_i$. Then $\pi_{LS}$ is a group model for $\mathcal{F}$. 
\end{Theorem}
%
The group models of Robinson type are iterated amalgams of
automorphism groups in the linking system over the $S-$normalizers
of the respective $\mathcal{F}$-centric subgroups of $S$.
\begin{Theorem}[\cite{Robinson1}, Theorem 2]
Let $\mathcal{F}$ be an Alperin fusion system on a finite $p-$group $S$ with $Aut_{\mathcal{L}}(P_1),...,Aut_{\mathcal{L}}(P_n)$.
The group
$\Gamma_R=Aut_{\mathcal{L}}(P_1)\underset{N_S(P_2)}{*}Aut_{\mathcal{L}}(P_2)\underset{N_S(P_3)}{*}\cdot.\underset{N_S(P_n)}*Aut_{\mathcal{L}}(P_1)$
is a group model for $\mathcal{F}$ where the maps used to define the amalgamations are $\delta_S: N_S(R_j)\rightarrow \Aut_\mathcal{L}(S)$ and $\delta_{R_j}: N_S(R_j)\rightarrow \Aut_\mathcal{L}(R_j)$.
\end{Theorem}
Corresponding to the various versions of Alperin's fusion theorem ($\mathcal{F}$-essential subgroups, $\mathcal{F}$-centric subgroups, $\mathcal{F}$-centric-radical subgroups) there exist several choices for the groups generating $\mathcal{F}$.\\

The group model constructed by the second author together with
Libman [\cite{Ian+Radu}, Theorem 1.1 and Proposition 4.1] is related to the normalizer decomposition and can be
described as follows.
\begin{Theorem}
Let $\mathcal{F}$ be a saturated fusion system on the finite $p-$group $S$. The finitely generated 
group 
$\mathcal{G}=Aut_{\mathcal{L}}(S)\underset{Aut_{\mathcal{L}}(P_2<S)}{*}Aut_{\mathcal{L}}(P_2)\cdots\underset{Aut_{\mathcal{L}}(P_n<S)}*Aut_{\mathcal{L}}(P_n)$
is a group model for $\mathcal{F}$ where $Aut_{\mathcal{L}}(P_i<S)=N_{Aut_{\mathcal{L}}(S)}(P_i)$ for $i=1,\cdots ,n$ and we amalgamate via the natural inclusion maps.
\end{Theorem}
The finite group $G=S\wr
\Sigma _{e(X)}$ constructed by Park \cite{Sejong} is not a group model since $S$ in the general case is not a Sylow $p-$subgroup of $G$, where $e(X)$ is a certain characteristic biset associated with $\mathcal{F}$.\\

We define a group model to be \textbf{minimal} if the map in the following theorem is an isomorphism. 
\begin{Theorem}[\cite{gmffs}]
Let $\mathcal{F}$ be a saturated fusion system over the finite
$p-$group $S$ and $\mathcal{G}$ a group model for $\mathcal{F}$ .
Then there exist a natural map of unstable algebras
$H^*(B\mathcal{G})\overset{q}{\rightarrow} H^*(\mathcal{F})$ making
$H^*(\mathcal{F})$ a module over $H^*(B\mathcal{G})$.
\end{Theorem}
We have a formula for the cohomology ring $H^*(B\mathcal{G};\mathbb{F}_p)$ which is always $F$-isomorphic to $H^*(\mathcal{F})$.
\begin{Theorem}[\cite{assclassfisoq}, Theorem 3.6] Let $\mathcal{F}$ be a saturated fusion system over the finite $p-$group 
$S$ and $\mathcal{G}$ one of the above group models for $\mathcal{F}$. Then there exist natural maps of algebras over the Steenrod
algebra $q: H^* (B\mathcal{G}) \rightarrow  H^* (\mathcal{F})$ and $r^*: H^*(\mathcal{F}) \rightarrow H^*(B\mathcal{G})$ such that we obtain a split
short exact sequence of unstable modules over the Steenrod algebra
$0\rightarrow W \rightarrow H^* (B\mathcal{G}) \overset{\leftarrow}{\rightarrow} H^* (\mathcal{F}) \rightarrow 0,$
where 
$W \cong Ker(Res^{\mathcal{G}}_S )\in\mathcal{N}il_1$ and the map $q$ is an $F$-isomorphism in the sense of Quillen \cite{Quillen}.
\end{Theorem}
The question of the existence of minimal group models for saturated fusion systems is related to the existence of the centric linking systems in the following way. One can define fusion systems and centric linking systems in a topological setting. We will need this when we make use of the fact that a group realizes a given fusion system if and only if its classifying space has a certain homotopy type.
 Recall \cite{BLO4} that an associated centric linking system exists if there exists a space which has the homology of the linking system. This motivates our search for minimal group models.\\ A \textbf{signalizer functor} in the sense of Aschbacher-Chermak on a group model $\mathcal{G}$ is an assignment $P\mapsto\theta (P)$ for every $\mathcal{F}-$centric subgroup $P\leq S$ such that $\theta (P)$ is a complement of $Z(P)$ in $C_{\mathcal{G}}(P)$ and such that if $gPg^{-1}\leq Q$ for $g\in\mathcal{G}$ then $\theta (Q)\leq g\theta (P)g^{-1}$. A signalizer functor gives rise to a centric linking system if it exists \cite{LS}. In \cite{LS} Libman and the author show that all previously discussed group models have a signalizer functor however without an explicit construction. This is one of the main results of this article.\\
All the group models presented above and the new construction in this article are equipped with a map $BG\rightarrow |\mathcal{L}|$ which is not necessarily true in the case in general. For a finite group $G$ we usually do not have a map $BG\rightarrow |\mathcal{L}_S^c(G)|$, the symmetric group on five elements is a counterexample for the classifying space of its $2$-local finite group. In the case that we have for a group $G$ and and moreover there is a map $BG\rightarrow |\mathcal{L}|$ which is an equivalence we solve the eighth problem of Oliver's list \cite{AKO} with  $G\cong\pi _1(|\mathcal{L}|)$. For this to happen the group $G$ does not have to be a group model for $\mathcal{F}$ as the linking system can be simply connected (e.g. the Solomon groups) but can be as for the general linear groups. There is no obvious pattern and we hope that our work on group models will allow to get better insight.
\subsection{Graphs of groups}
\label{cAppendix} We give a brief introduction to graphs of groups stating results we need.
A \textbf{finite directed graph} $\Gamma$
consists of two sets, the \textbf{vertices} $V$ and the
\textbf{directed edges} $E$, together with two functions $\iota,
\tau : E\rightarrow V$. For $e \in E,\iota (e)$ is called the
\textbf{initial vertex} of $e$ and $\tau (e)$ is the
\textbf{terminal vertex} of $e$. Multiple edges and loops are
allowed in this definition. The graph $\Gamma$ is \textbf{connected}
if the only equivalence relation on $V$ that contains all $(\iota
(e),\tau (e))$ is the relation with just one class. A graph $\Gamma$
may be viewed as a category, with objects the disjoint union of $V$
and $E$ and two non-identity morphisms with domain $e$ for each
$e\in E$, one morphism $e\rightarrow \iota (e)$ and one morphism
$e\rightarrow \tau (e)$. A \textbf{graph $\Gamma$ of groups} is a
connected graph $\Gamma$ together with groups $G_v, G_e$ for each
vertex and edge and injective group homomorphism
$f_{e,\iota}:G_e\rightarrow G_{\iota (e)}$ and $f_{e,\tau
(e)}:G_e\rightarrow G_{\tau (e)}$ for each edge $e$. 
Let $(S,\mathcal{F},\mathcal{L})$ be a $p$-local finite group, $P_1,\dots,P_n$ a collection of $\mathcal{F}$-centric subgroups of $S$ which might contain multiples, with a generating collection $K=\{K_1,\dots,K_n\}$ where $K_i\leq\Aut_\mathcal{L}(P_i)$. We define a \textbf{graph of groups $\Gamma_K$ associated with $K$}. 
%
%
Let $F_{i,j}$ be the subgroup of $\Aut_\mathcal{F}(\langle P_i,P_j\rangle)$ consisting of automorphisms that restrict in $\mathcal{F}$ to both an automorphism of $P_i$ contained in $\pi(K_i)$ and an automorphism of $P_j$ contained in $\pi(K_j)$. Define $\overline{K_{i,j}}$ as the preimage of $F_{i,j}$ in $\Aut_\mathcal{L}(\langle P_i,P_j\rangle)$. Note that $\overline{K_{i,j}}=\overline{K_{j,i}}$. 
By \cite[Prop 2.11.]{Lib06} we can define an injective restriction homomorphism $k_{i,j}:\overline{ K_{i,j}}\rightarrow K_j$.
 A graph of groups $\Gamma_K$ will be called a graph \textbf{associated} with $K$, if there exist Sylow $p$-subgroups $H_i\leq K_i$ and $H_{i,j}<K_{i,j}$ such that $\Gamma_K$ has $n$ vertices, with corresponding groups $K_i$ and such that there is a directed edge from $K_i$ to $K_j$ in $\Gamma_K$ if and only if $H_j<k_{i,j}(H_{i,j})$, and to this edge correspond the injective group homomorphisms $k_{j,i}:  K_{i,j}\rightarrow K_i$ and $k_{i,j}:  K_{i,j}\rightarrow K_j$ for a subgroup $H_{i,j}\leq K_{i,j}\leq \overline{ K_{i,j}}$.
In our applications we will need to choose a special subgraph of $\Gamma_K$.
An \textbf{aborescence} is a directed graph in which, for a vertex $u$ called the \textbf{root} and any other vertex $v$, there is exactly one directed path from $u$ to $v$.
 Let $\Gamma_K$ be a graph associated with $K$. A subgraph $T$ of $\Gamma_K$ which is an aborescence with root $t$ is called a \textbf{generating tree of $\Gamma_K$}, if $P_t\vartriangleleft S$ and $H_t=S<K_t$.
It might not always be possible to find a generating tree for $\Gamma_K$.
\section{A new family realizing saturated fusion systems}
We aim to give a family of group models for saturated fusion systems. This family of group models contains versions of classical group models as special cases. The main advantage over previous constructions will be discussed at the end of this article. 
\begin{Theorem}
\label{groupmodelkt}
 Let $K=\{K_1,\dots,K_n\}$ be a generating collection, $\Gamma_K$ the associated graph for some choice of $H$'s and $T$ a generating tree in $\Gamma_K$. Then the amalgam $\pi_{K,T}$ over the graph of groups $T$ is a group with the following properties:
\begin{enumerate}
 \item $S$ is a Sylow $p$-subgroup of $\pi_{K,T}$.
 \item $\mathcal{F}=\mathcal{F}_S(\pi_{K,T})$.
 \item There is a signalizer functor on $\pi_{K,T}$ which induces $\mathcal{L}$.
\item $H^*(|\mathcal{L}|,\F_p)$ is a retract of $H^*(B\pi_{K,T})$ in the category of unstable algebras. It is equal to the image of $H^*(\pi_{K,T},\F_p)\rightarrow H^*(S,\F_p)$, the product of any two elements
 in the kernel is zero.
 \item $B\pi_{K,T}$ is $p$-good.
 \item $H^*(B\pi_{K,T})$ is finitely generated.
 \item $|\mathcal{L}|\pcom$ is a stable retract of $(B\pi_{K,T} )\pcom $.
 \item $H^*(B\pi_{K,T})$ is $F-$isomorphic to the stable elements in
 the sense of Quillen.
\end{enumerate}
\end{Theorem}
\underline{Proof:} The proof is based on results of Libman and ourselves \cite{LS}. Consider the category $\mathcal{C}$
consisting of objects $l_0,\dots,l_n$ as well as $l_{i,j}$ for those
$i,j$ such that $K_{i,j}$ is in $T$. Those objects correspond to the
groups mentioned above (i.e. $K_0,\dots K_n$ and the  $K_{i,j}$ in
$T$). In $\mathcal{C}$ let there be a unique morphism from
$l_{i,j}$ to $l_i$ and $l_j$. We have a functor $\gamma$ from this
category into $\hoTop$, which sends $l_i$ to $BK_i$ and $l_{i,j}$ to
$BK_{i,j}$, and the morphisms to those induced by the monomorphisms
$k_{j,i}$ and $k_{i,j}$. We can include each of the classifying spaces into
$|\mathcal{L}|$ (induced by an inclusion of categories
$\mathcal{B}K_i\rightarrow
\mathcal{B}\Aut_\mathcal{L}(P_i)\rightarrow  \mathcal{L}$). We now
want to show that this commutes up to homotopy with the the
morphisms from $\mathcal{C}$, i.e. we want the outer diagram to
commute up to homotopy:
\[\begin{xy}
   \xymatrix{
     & BK_i\ar[rd] &  \\
    BK_{i,j}\ar[ru]\ar[rd]\ar[r] & B\Aut_\mathcal{L}(\langle P_i,P_j \rangle)\ar[r] & |\mathcal{L}| .\\
     & BK_j\ar[ru] &
  }
\label{hocolim_diagram}
 \end{xy}\]

To see this consider the functors
$F_1: \mathcal{B}K_{i,j}\rightarrow\mathcal{B}K_i\rightarrow
\mathcal{L}$ and $F_2:  \mathcal{B}K_{i,j}\rightarrow
\mathcal{B}\Aut_\mathcal{L}(\langle P_i,P_j \rangle \rightarrow
\mathcal{L}$. There is a natural transformation $\eta:
F_1\rightarrow F_2$, which maps the unique object to the morphism
$\widehat{e}=\delta(e)\in \Mor_{\mathcal{L}}(P_i,\langle
P_i,P_j\rangle)$. 
The following diagram commutes by the definition of $k_{j,i}$:
\[\begin{xy}
   \xymatrix{
P_i \ar[r]^{\widehat{e}}\ar[d]_{k_{j,i}(f)=F_1(f)} & \langle P_i,P_j\rangle \ar[d]^{F_2(f)=f}\\
P_i\ar[r]_{\widehat{e}} \ar[r] &\langle P_i,P_j\rangle .
  }
 \end{xy}\]
Thus $\eta$ is a natural transformation. Taking realizations shows that the two triangles in the diagram above are homotopy commutative.
By the universal property, this gives a map
$f: \hocolim \gamma \rightarrow |\mathcal{L}|$. We claim that the composition with
$|\mathcal{L}|\rightarrow |\mathcal{L}|_p^\wedge$ and restriction to
$BS\leq BK_t$ is homotopic to $\theta: BS\rightarrow
|\mathcal{L}|_p^\wedge$ (where we use notation as in \cite{LS}). To see this, consider the following diagram
 \[\begin{xy}
     \xymatrix{
      & & & BK_t\ar[rd] &  \\
     {B}S\ar[rr]^{{B}\delta_S}\ar[rrru]^{{B}\delta_{P_t}}\ar[rrrd]_{{B}\delta_S} & &{B}K_{t,0}\ar[ru]_{{B}k_{0,t}}\ar[rd]^{{B}k_{t,0}}& & |\mathcal{L}|\\
      & & & {B}\Aut_\mathcal{L}(S) \ar[ru] &
    }
  \end{xy}\]
where we set $K_0=\Aut_\mathcal{L}(S)$ for the purpose of this diagram. The morphism $\delta_{P_t}: S\rightarrow \Aut_\mathcal{L}(P_t)$ makes sense because we required $P_t\vartriangleleft S$. The image lies in $K_t$ by another condition on generating trees. Thus the upper left morphism makes sense. The diagram 
\[\begin{xy}
   \xymatrix{
P_t \ar[r]^{\widehat{e}}\ar[d]_{\delta_{P_t}(g)} & S \ar[d]^{\delta_S(g)}\\
P_t\ar[r]_{\widehat{e}} \ar[r] &S
  }
 \end{xy}\]
commutes for all $g\in S$ because $\delta$ is a functor and $\widehat{e}=\delta_{P_t,S}(e)$. Thus $k_{0,t}\circ \delta_S=\delta_{P_t}$, and so the upper triangle strictly commutes.
We have $\langle P_t,S\rangle =S$, so $K_{t,0}$ is a subgroup of $\Aut_\mathcal{L}(S)$, so the restriction $k_{t,0}$ is really the inclusion. Thus the lower triangle also commutes strictly, consisting of morphisms whose restriction to $P_t$ is in $K_t$. 
By the reasoning employed at the beginning of this proof, the right square commutes up to homotopy. Composition with $|\mathcal{L}|\rightarrow|\mathcal{L}|_p^\wedge$ yields the claim.
We have a graph of groups $T$ which we can regard as a category. Note that the functor $B\mathcal{G}: T\rightarrow \hoTop$ is
isomorphic to $\gamma$ discussed above. All the groups involved have
a Sylow $p$-subgroup and the map from a chosen Sylow $p$-subgroup of
$K_{i,j}$ to a Sylow $p$-subgroup of $K_j$ is surjective if $K_{i,j}$ is in $T$. Furthermore, $T$ is a tree and contains a path from $v_0=K_t$ to every other vertex. Thus we can apply \cite[Proposition 3.3]{LS}.
This immediately proves claims 1, 5 and the last part of 4, as well as the fact that  $\hocolim \gamma \simeq B\pi_{K,T}$.
By combining this with the result above, we get a map $f:B\pi_{K,T}\simeq \hocolim \gamma\rightarrow  |\mathcal{L}|_p^\wedge$, whose restriction to $BS$ is homotopic to $\theta$.
We want to use \cite[Theorem 1.1]{LS}, so we need $\mathcal{F}\subset \mathcal{F}_S(\pi_{K,T})$. This is the case, as $\pi_{K,T}$ contains all $K_i$, so $\mathcal{F}_S(\pi)$ contains all $K_i'$. But then $\mathcal{F}_S(\pi_{K,T})$ also contains the fusion system generated by the $K_i'$, which is $\mathcal{F}$, as $K$ is generating. Thus \cite[Theorem 1.1]{LS} yields that the group model $\pi_{K,T}$ also has the properties 2,3 and 4. 
The group $\pi _{K,T}$ is a finite amalgam of finite groups which is generated by elements of $p'$-order and $S$. Let $M$ be the subgroup of $\pi _{K,T}$ generated by all elements of $p'$ -order.
Note that $\pi _{K,T}$ and $S$ surjects on $\pi _{K,T} / M$ and therefore $\pi _{K,T} / M$ is a finite $p$-group. The group $M$ is $p$-perfect since it is generated
by $p'$-elements. Let $X$ be the cover of $B\pi _{K,T}$  with fundamental group $M$. Using the results from
\cite[ VII.3.2]{BK}, we have that $X$ is $p$-good and $X\pcom$
is simply connected. Hence the sequence $X\pcom\rightarrow (B\pi _{K,T})\pcom\rightarrow B(\pi _{K,T} / M)$ is a fibration sequence
and so $(B\pi _{K,T})\pcom$ is $p$-complete by \cite[II.5.2(iv)]{BK}. So $B\pi _{K,T}$ is $p$-good. 
 Recall that
   $H_i\in Syl_p(K_i)$ for all $i=1,...,n$. It follows from \cite[Lemma 2.3.]{BLO1} and
    \cite[Theorem 4.4.(a)]{BLO2} that $H^*(B(K_i))$ is finitely generated over $H^*(|\mathcal{L}|)$ for
    all $i=1,...,n$, and $H^*(|\mathcal{L}|)$ is noetherian as follows from
     \cite[Proposition 1.1. and Theorem 5.8.]{BLO2}. Therefore the Bousfield-Kan spectral sequence
     for $H^*(BG)$ is a spectral sequence of finitely generated $H^*(|\mathcal{L}|)-$modules,
     the $E_2$ term with $E_2^{s,t}= \underset{\mathcal{C}}{lim^s}H^t(F(-);\mathbb{F}_p)$ is concentrated
     in the first two columns and $E_{2}=E_{\infty}$ for placement reasons. Therefore $H^*(B\pi_{K,T})$ is a
      finitely generated module over $H^*(|\mathcal{L}|)$. 
For property 7 recall we have a commutative diagram
\xymatrix@R=9pt@C=9pt{&{\Sigma ^{\infty}BS\pcom}\ar[1,1]^{\Sigma ^{\infty}Bincl\pcom}\ar[1,-1]_{\Sigma ^{\infty}B(\delta _S)\pcom}&{}\\
{\Sigma ^{\infty}|\mathcal{L}|\pcom}&&{\Sigma ^{\infty}(B\pi_{K,T})\pcom.}\ar[0,-2]^{\Sigma ^{\infty}q\pcom}\\
} By Ragnarsson's work  \cite{Ragnarsson} there is a map of spectra $\sigma _{\mathcal{F}}:\Sigma ^{\infty}|\mathcal{L}|\pcom\rightarrow\Sigma ^{\infty}BS \pcom$ such that the composition $ \Sigma ^{\infty}|\mathcal{L}|\pcom\overset{\sigma _{\mathcal{F}}}{\longrightarrow}\Sigma ^{\infty}BS\pcom\overset{\Sigma ^{\infty}B(\delta _S)\pcom}{\longrightarrow}\Sigma ^{\infty}|\mathcal{L}|\pcom$ is the identity. 
Since $\Sigma ^{\infty}B(\delta _S)\pcom\circ \sigma _{\mathcal{F}}
=\Sigma ^{\infty}q\pcom\circ\Sigma ^{\infty}Bincl\pcom\circ \sigma
_{\mathcal{F}}$ we have that the simplicial set $|\mathcal{L}|\pcom$ is a stable retract of
$(B\pi_{K,T})\pcom$. 
Point 8 follows from the definition of $F-$isomorphism and 4.$\Box$\\ 

We can use group models to compute the fundamental group of $|\mathcal{L}|\pcom$ under certain conditions.
\begin{Proposition}
\label{BG SIMEQ |L|}
If $f:B\pi_{K,T}\rightarrow  |\mathcal{L}|_p^\wedge$ is a homotopy equivalence then $\pi _1(|\mathcal{L}|_p^\wedge)\cong\pi _{K,T} / M $.
\end{Proposition}
\underline{Proof:} The statement follows from \cite[Prop 1.11]{AKO} since in this case $M$ is the maximal $p$-perfect subgroup of $\pi _{K,T}\cong\pi _1(|\mathcal{L}|)$. $\Box$

\begin{Conjecture}
Let $\mathcal{F}$ be a saturated fusion system over the finite $p$-group $S$ and $\mathcal{G}$ a group model for $\mathcal{F}$ with trivial signalizer functor ( which implies $\mathcal{G}$ contains subgroups isomorphic to $Aut_{\mathcal{L}}(P)$ for all $\mathcal{F}$-centric subgroups $P\leq S$). Then $\pi_1 (|\mathcal{L}|\pcom)\cong B\mathcal{G}\pcom$.
\end{Conjecture}
\subsection{The group model $\Gamma_R$}
Let $(S,\mathcal{F},\mathcal{L})$ be a $p$-local finite group. Let $R_1,\dots, R_n$ be fully normalized representatives of isomorphism classes of the centric-radical subgroups of $S$. Fix $R_1=S$, set $K_i=H_i=N_S(R_i)$. In this case $\Gamma_K$ contains a directed star graph $T$ with root $\Aut_\mathcal{L}(S)$. We obtain a group model $\Gamma _R\cong\Aut_\mathcal{L}(S)\underset{N_S(R_i)}{*}\Aut_\mathcal{L}(R_1)\ast\dots\underset{N_S(R_i)}{*}\Aut_\mathcal{L}(R_n)\cong\pi_R$
, recovering results from Robinson
 \cite{Robinson1}.
Theorem \ref{groupmodelkt} shows that it was not
 crucial that $\Aut_\mathcal{L}(S)$ was used as root - any $K_i$ works as long as $R_i\vartriangleleft S$.
\subsection{The group model $\pi_R$}
Let $(S,\mathcal{F},\mathcal{L})$ be a $p$-local finite group, $R_1,\dots, R_n$ $\mathcal{F}$-centric subgroups of $S$. Fix $R_1=S$, set $K_i=\Aut_\mathcal{L}(R_i)$, $K_{i,j}=\overline{K_{i,j}}$ and $H_i=N_S(R_i)$. In this case $\Gamma_K$ contains a directed star graph $T$ with root $\Aut_\mathcal{L}(S)$. We obtain  $\pi_{K,T}\cong\Aut_\mathcal{L}(S)\underset{\Aut_\mathcal{L}(R_1<S)}{*}\Aut_\mathcal{L}(R_1)\ast\dots\underset{\Aut_\mathcal{L}(R_n<S)}{*}\Aut_\mathcal{L}(R_n)\cong\pi_R$
which is a group model for $\mathcal{F}$ recovering results from Libman and the second author
 \cite{LS}.
Theorem \ref{groupmodelkt} shows that it was not
 crucial that $\Aut_\mathcal{L}(S)$ was used as root - any $K_i$ works as long as $R_i\vartriangleleft S$.
\subsection{Splitting up  the $\Aut_{\mathcal{L}}(R_i)$}
We can generalize $\pi_R$ in the following way. Let $R_1,\dots,R_n$ be as above, and $L_i=\Aut_\mathcal{L}(R_i)$. Assume $K_1,\dots,K_m$ to be a collection of subgroups of the $L_i$'s (ordered so that $K_{j_i},\dots, K_{j_{i+1}-1} < L_i$ for suitable $j$'s, with $j_1=1, j_{n+1}=m+1$), and such that $L_i$ is generated by $K_{j_i},\dots, K_{j_{i+1}-1}$ and $N_S(R_i)<K_j$ for $j_i \leq j < j_{i+1}$.
In this case, $K$ is a generating collection.
Let $P_j=R_i$ if $K_j<L_i$. Clearly $P_1=S$. Thus $K_{1,j}$
consists of the elements of $K_1$ whose restriction to $P_j$ is in $K_j$, or in
other words the elements of $K_j$, which are extendable to
automorphisms of $S$ that lie in $K_1$. By the third axiom of linking
systems, all elements of $N_S(P_j)<\Aut_\mathcal{L}(P_j)$ are extendable, and their extension lies in $S<\Aut_\mathcal{L}(S)$, so in $K_1$. We can chose  $H_j=N_S(P_j)$, and our assumptions as well as
the comment just made show that $\Gamma_K$ contains a directed star graph $T$ with root $K_1$ and edges directed away from it. As $P_1=S$, $T$ is a generating tree, and we get a group model:
\[\pi_{K,T}=K_1\underset{K_{1,2}}{*}K_2\ast\dots\underset{K_{1,m}}{*} K_m.\]
\subsection{The group model $\pi_{LS}'$}
%
%
We can modify the Leary-Stancu group model in the following way.
\begin{Theorem}
 Let $(S,\mathcal{F},\mathcal{L})$ be a $p$-local finite group, and let $\phi_1,\dots,\phi_n$ generate $\mathcal{F}$. Assume as we can that all $\phi_i$ have order coprime to $p$ and that they are automorphisms of $\mathcal{F}$-centric subgroups of $S$. Then $\pi_{LS}'=S\ast Fr(t_1,\dots, t_n)/\langle t_i^{\ord(\phi_i)}=1, t_i u t_i^{-1}=\phi_i(u)\rangle$ is a group model of the type described in Theorem \ref{groupmodelkt}.
\end{Theorem}
\underline{Proof:}
Denote the subgroup of $S$ that $\phi_i$ is an automorphism of by $P_i$.
Let $K_i'=\langle P_i/Z(P_i),\phi_i\rangle$ as a subgroup of
$\Aut_\mathcal{F}(P_i)$, and let $K_i$ be the preimage in
$\Aut_\mathcal{L}(P_i)$. 
The kernel of $\pi:\Aut_\mathcal{L}(P_i)\rightarrow \Aut_\mathcal{F}(P_i)$ is a $p$-group, so as $\phi_i$ has order prime to $p$ we can find a preimage $t_i$ of $\phi_i$ which also has order $\ord(\phi_i)$. Clearly $K_i=\langle P_i, t_i\rangle$. Set $K_0=S<\Aut_\mathcal{L}(S)$. By construction, the collection $K_0,\dots K_n$ is generating. By the third axiom of linking spaces, $t_iut_i^{-1}=\phi_i(u)$ for $u\in P_i$. In particular, $P_i\vartriangleleft K_i$. Thus the order of $K_i$ is $|P_i|\cdot \ord(\phi_i)$, so we can set $H_i=P_i$ and $H_0=S$. Clearly $K_{0,i}=P_i$, so the directed star graph $T$ with root $K_0$ is a generating tree. Theorem \ref{groupmodelkt} yields a group model $\pi_{K,T}=S\underset{P_1}{
*} \langle P_1,
t_1\rangle \ast \dots \underset{P_n}{*} \langle P_n , t_n\rangle$. 
We claim that this group is $\pi_{LS}'$ as given in the statement of this proposition.
Consider the obvious group homomorphism $S\ast Fr(t_1,\dots, t_n) \rightarrow \pi_{K,T}$. It was already noted that the relations which are divided out for $\pi_{LS}'$ also hold in $\pi_{K,T}$. Thus this homomorphism descends to $\psi: \pi_{LS}'\rightarrow \pi_{K,T}$. For the inverse direction, first note that $K_i=\langle P_i, t_i\rangle \cong P_i \ast C_{\ord(\phi_i)} \langle t_i\rangle/\langle t_iut_i^{-1}=\phi_i(u) \rangle$, as there is a surjective group homomorphism from the right side to the left and both groups have the same order. Thus the obvious group homomorphism $K_i \rightarrow \pi_{LS}'$ is well-defined. Clearly all these homomorphisms make the necessary diagrams commute so that by the universal property of pushouts they combine to a group homomorphism $\psi^{-1}: \pi_{K,T}\rightarrow \pi_{LS}'$ which is an inverse for $\psi$. Thus $\pi_{K,T}\cong \pi_{LS}'$. $\Box$
\section{Signalizer functor constructions}
Each $p$-local finite group $(S,\mathcal{F},\mathcal{L})$ can be realized by some group and some choice of signalizer functor in the sense of Aschbacher-Chermak. This was proved in \cite{LS} however without an explicit construction which we give in this section.
\subsection{Signalizer functor for $\Gamma _R$}

Let $(S,\mathcal{F},\mathcal{L})$ be a $p$-local finite group and $R_1,\cdots , R_n$ be $\mathcal{F}$-centric subgroups of $S$ such that the group $\Gamma_R=\Aut_{\mathcal{L}}(S)\ast\dots\underset{N_S(R_n)}{\ast}\Aut_{\mathcal{L}}(R_n)$ is a Robinson group model for $\mathcal{F}$.
We construct a signalizer functor for $\Gamma_R$ using the algebraic structure of $\mathcal{L}$. This means choosing subgoups $C^*_{\Gamma _R}(P)\leq C_{\Gamma _R}(P)$ for each $P\in\mathcal{F}^c$ such that $Mor_{\mathcal{L}}(P,Q)\cong N_{\Gamma _R}(P,Q)/C^*_{\Gamma _R}(P)$ for each $P$ and $Q$ which are $\mathcal{F}-$centric subgroups of $S$ and provides the complete solution to Oliver's seventh problem \cite{AKO} whether every centric linking system is the centric linking system of a group. The signalizer functor is defined as the kernel of a collection of maps $\{ \psi_P: Aut_{\tau^c_S(\Gamma_R)}(P)\rightarrow \Aut_\mathcal{L}(P) \}_{P\in \mathcal{F}^c}$. 
 The difference to prior results obtained by Libman and the author \cite{LS} is that they define what we call $\{ \psi_P: Aut_{\mathcal{T}^c_S(\Gamma_R)}(P)\rightarrow \Aut_\mathcal{L}(P) \}_{P\in \mathcal{F}^c}$ via topological maps, whereas we use the algebraic structure of the categories involved directly.\\
 We write $\delta_{S}(N_S(R_j))$ or
 $\delta_{R_j}(N_S(R_j))$ instead of $N_S(R_j)$ to distinguish
 between the elements $N_S(R_j)$ as a part of $\Aut_{\mathcal{L}}(S)$
 and $\Aut_{\mathcal{L}}(R_j)$.
The maps used to define the amalgamations are $\delta_S: N_S(R_j)\rightarrow \Aut_\mathcal{L}(S)$ and $\delta_{R_j}: N_S(R_j)\rightarrow \Aut_\mathcal{L}(R_j)$.
We define $\psi_P: \Aut_{\mathcal{T}^c_S(\Gamma_R)}(P)\rightarrow \Aut_\mathcal{L}(P)$ for every $\mathcal{F}$-centric $P\leq S$.
Here, $\Aut_{\mathcal{T}^c_S(\Gamma_R)}(P)=N_{\Gamma_R}(P)$, where $P$ is embedded into $\Gamma_R$ via $P\leq S \leq \Aut_\mathcal{L}(S)\leq \Gamma_R$.
From \cite[Lemma 1]{Robinson1} applied to iterated amalgams it follows that elements of $N_{\Gamma_R}(P)$ must have the form $a_1\dots a_k$ such that there exists subgroups $P_i\leq S$,  $i=1,\cdots k$ with:
$a_i\in \Aut_\mathcal{L}(R_{j(i)})$ for a $j(i)$,
$P_i,P_{i-1}\leq N_S(R_{j(i)})$,  
$P_k=P$,
$a_i\delta_{R_{j(i)}}({P_{i}})a_i^{-1} = \delta_{R_{j(i)}}(P_{i-1})$,
$P_0=P$.
\begin{Proposition}
 Let $g=a_1a_2\dots a_k$ be such an element in $N_{\Gamma_R}(P)$ as above. Then for every $i$, $a_i$ has a unique extension $\widetilde{a_i}: \langle R_{j(i)},P_i\rangle \rightarrow \langle R_{j(i)}, P_{i-1}\rangle $ as a morphism in $\mathcal{L}$. Furthermore, $\widetilde{a_i}$ is an isomorphism and has a unique restriction $\widehat{a_i}: P_i\rightarrow P_{i-1}$, which is also an isomorphism.
\end{Proposition}
\underline{Proof:} We have $a_i\in N_{\Gamma_R}(P_i,P_{i-1})$ and $a_i\in \Aut_\mathcal{L}(R_{j(i)})$. As $P_i\leq N_S(R_{j(i)})$ and $P_{i-1}\leq N_S(R_{j(i)})$ we get $R_{j(i)}\vartriangleleft \langle R_{j(i)}, P_i\rangle$ and $R_{j(i)}\vartriangleleft \langle R_{j(i)}, P_{i-1}\rangle$. Furthermore as the way we identify $N_S(R_{j(i)})$ in $\Aut_{\mathcal{L}}(S)$ and $\Aut_{\mathcal{L}}(R_{j(i)})$ is via $\delta_{S,S}$ and $\delta_{R_{j(i)},R_{j(i)}}$, we have:
$a_i\delta_{R_{j(i)}}(P_i)a_i^{-1}=\delta_{R_{j(i)}}(P_{i-1})\leq\delta_{R_{j(i)}}(\langle R_{j(i)},P_{i-1}\rangle)$ and
$a_i\delta_{R_{j(i)}}(R_{j(i)})a_i^{-1}\leq \delta_{R_{j(i)}}(R_{j(i)})\leq\delta_{R_{j(i)}}(\langle R_{j(i)},P_{i-1}\rangle)$, where the first part of the second inequality follows from the third axiom of linking systems.
We obtain $a_i\delta_{R_{j(i)}}(\langle R_{j(i)},P_i\rangle)a_i^{-1}\leq  \delta_{R_{j(i)}}(\langle R_{j(i)},P_{i-1}\rangle)$
and therefore \cite[Proposition 4 (c)]{O4} can find a unique extension $\widetilde{a_i}$ in the diagram in $\mathcal{L}$ below:
\[\begin{xy}
   \xymatrix{
R_{j(i)} \ar[r]^{a_i}\ar[d]_{\iota} & R_{j(i)} \ar[d]^{\iota}\\
\langle R_{j(i)},P_{i}\rangle\ar[r]_{\widetilde{a_i}} \ar[r] &\langle R_{j(i)},P_{i-1}\rangle .
  }
 \end{xy}\]
Analogously we obtain $\widetilde{a_i^{-1}}$ for $a_i^{-1}$ which has the property
$a_i^{-1}\delta_{R_{j(i)}}(P_{i-1})a_i=\delta_{R_{j(i)}}(P_i)$.
 After extending $\delta_{R_{j(i)}}(1)$ and uniqueness we obtain that $\widetilde{a_i}$ is an isomorphism in $\mathcal{L}$ with
  inverse $\widetilde{a_i^{-1}}$. 
%
 In order to define $\widehat{a_i}:P_i\rightarrow P_{i-1}$ as the unique restriction of $\widetilde{a_i}$, by \cite[Proposition 4 (b)]{O4} we have to check that $\pi(\widetilde{a_i})(P_i)\leq P_{i-1}$.
  So let $p\in P_i$. Because $a_i \delta_{R_{j(i)}}(P_i)a_i^{-1}=\delta_{R_{j(i)}}(P_{i-1})$,
  there is a $p'\in P_{i-1}$ with $a_i \delta_{R_{j(i)}}(p)a_i^{-1}=\delta_{R_{j(i)}}(p')$. 
   Consider the
  following diagram in $\mathcal{L}$:
%

\[
\begin{tikzcd}[row sep=scriptsize, column sep=scriptsize]
& R_{j(i)} \arrow{ddl}[swap]{\delta_{R_{j(i)}}(p)}\arrow{rr}{a_i}\arrow{ddd} & & R_{j(i)} \arrow{ddl}[swap, near end]{\delta_{R_{j(i)}}(p')}\arrow{ddd} \\
\\
R_{j(i)} \arrow[crossing over]{rr}[near end]{a_i}\arrow{ddd} & & R_{j(i)} \\
& \langle P_i,R_{j(i)}\rangle \arrow{ddl}[near start]{\delta_{\langle R_{j(i)},P_i\rangle}(p)}\arrow{rr}[near start]{\widetilde{a_i}} & & \langle P_{i-1},R_{j(i)}\rangle \arrow{ddl}[near start]{\delta_{\langle R_{j(i)},P_{i-1}\rangle}(p')} \\
\\
\langle P_i,R_{j(i)}\rangle \arrow{rr}[swap]{\widetilde{a_i}} & & \langle P_{i-1},R_{j(i)}\rangle \arrow[crossing over, leftarrow]{uuu}\\
\end{tikzcd}
\]

where the four unlabeled morphisms are $\iota=\delta(1)$.
The top square commutes by definition of $p'$. The squares on the left and right commute because $\delta$
is a functor. The front and back squares commute by definition of $\widetilde{a_i}$. In $\mathcal{L}$
all morphisms are epimorphisms \cite[Proposition 4 (d)]{O4}, making extensions unique, so
$\delta_{\langle R_{j(i)},P_{i-1}\rangle}(p')\circ \widetilde{a_i}$ is an extension of
 $\delta_{R_{j(i)}}(p')\circ a_i$, which is the same as $ a_i\circ\delta_{R_{j(i)}}(p)$, which has extension
  $\widetilde{a_i}\circ\delta_{\langle R_{j(i)},P_{i}\rangle}(p) $. Thus the lower square also commutes.
We want to show that $\pi(\widetilde{a_i})(p)=p'$. As $\delta$ is
injective on morphism sets, this follows now from $\delta_{\langle
R_{j(i)},P_{i-1}\rangle}(\pi(\widetilde{a_i})(p))=\widetilde{a_i}\circ
\delta_{\langle R_{j(i)},P_{i}\rangle}(p)\circ
\widetilde{a_i}^{-1}=\delta_{\langle R_{j,i},P_{i-1}\rangle}(p')$,
where in the first step the third axiom of linking systems was used.
Thus the restriction $\widehat{a_i}$ of $\widetilde{a_i}$ exists and
is unique. 
As above we can get a restriction $\widehat{a_i^{-1}}$ of $\widetilde{a_i}^{-1}$, and by uniqueness of the restriction of $\widetilde{id}$ we get that $\widehat{a_i^{-1}}=\widehat{a_i}^{-1}$, so $\widehat{a_i}$ is an isomorphism in $\mathcal{L}$.
$\Box$
%
%
%

This allows to define a map $\psi_P: \Aut_{\mathcal{T}^c_S(\Gamma_R)}(P)\rightarrow \Aut_\mathcal{L}(P)$ for every $\mathcal{F}$-centric $P\leq S$.
%
%
%
%
\begin{Proposition}
For every $a_1a_2\dots a_n\in N_{\Gamma_R}(P)$ given in the form above set $\psi_P(a_1a_2\dots a_n)=\widehat{a_1}\circ\widehat{a_2}\circ\dots\circ\widehat{a_n}$ where on the right we compose morphisms in $\mathcal{L}$. This is well-defined and makes $\psi_P$ into a group homomorphism  $\psi_P: \Aut_{\mathcal{T}^c_S(\Gamma_R)}(P)\rightarrow \Aut_\mathcal{L}(P)$.
\end{Proposition}
\underline{Proof:} Passing from $a_1,\dots,a_n$ to $\psi_P(a_1\dots a_n)$ is well-defined. We only need to show independence of the factorization of $g\in N_{\Gamma_R}(P)$. By Serre's form of reduced words \cite[Theorem 1]{Serre} for this it 
suffices to show that for $a\in N_S(R_i)\cap N_S(R_j)$ we have that $\widehat{a}: Q\rightarrow Q'$ is independent of the choice of $\delta_{R_i}(a)\in\Aut_\mathcal{L}(R_i)$ or $\delta_{R_j}(a)\in\Aut_\mathcal{L}(R_j)$, and that if $b,c\in\Aut_\mathcal{L}(R)$, then $\widehat{b\cdot c}=\widehat{b}\circ\widehat{c}$.
If $a\in N_S(R)$ and $\delta_R(a)\delta_R(Q)\delta_R(a)^{-1}=\delta_R(Q')$, so 
 $aQa^{-1}=Q'$, then $a\in N_S(\langle R,Q\rangle,\langle R,Q'\rangle)$. As $\delta$ is a functor, $\delta_{\langle R,Q\rangle,\langle R,Q'\rangle}(a)$ is an extension of  $\delta_{R}(a)$ and is thus the unique choice for $\widetilde{\delta_{R}(a)}$. Restricting $\delta_{\langle R,Q\rangle,\langle R,Q'\rangle}(a)$ to $Q,Q'$ we obtain  $\delta_{Q,Q'}(a)$. Thus $\widehat{\delta_R(a)}=\delta_{Q,Q'}(a)$. This does not depend on $R$, so yields the same result for $\delta_{R_i}(a)\in\Aut_\mathcal{L}(R_i)$ and  $\delta_{R_j}(a)\in\Aut_\mathcal{L}(R_j)$ in case that $a\in N_S(R_i)\cap N_S(R_j)$.
Let $b,c\in\Aut_\mathcal{L}(R)$, and $Q'',Q',Q \leq N_S(R)$, with
$cQc^{-1}=Q',bQ'b^{-1}=Q''$.
 By joining the two respective
commutative diagrams 
\[\begin{xy}
   \xymatrix{
R \ar[r]^{b}\ar[d] & R \ar[d]\ar[r]^{c} & R \ar[d]\\
\langle R,Q\rangle\ar[r]^{\widetilde{b}} & \langle R,Q'\rangle\ar[r]^{\widetilde{c}} & \langle R,Q''\rangle \\
Q \ar[u]\ar[r]_{\widehat{b}}& Q' \ar[u]\ar[r]_{\widehat{c}}& Q''\ar[u]
  }
 \end{xy}\]
 
 we can easily see that
$\widetilde{b}\circ\widetilde{c}$ is an extension of $b\circ c$, so
by uniqueness $\widetilde{bc}=\widetilde{b}\circ\widetilde{c}$. Restriction of
$\widetilde{b}\circ\widetilde{c}$ is the composition of the
individual restrictions, 
so the statement follows.
Since
$\psi(a_1\cdot\dots \cdot a_n\cdot b_1\cdot\dots \cdot c_m)=\widehat{a_1}\circ\dots\circ\widehat{a_n}\circ\widehat{b_1}\circ\dots\circ\widehat{b_m}=\psi(a_1\dots a_n)\circ\psi(b_1\dots b_m)$
we have:
 $\psi_P$ is a group homomorphism. $\Box$
\begin{Proposition}
\label{signalizer_functor}
 Let $P\leq S$ be $\mathcal{F}$-centric and let $\Theta(P)=\Ker\left(\psi_P: \Aut_{\mathcal{T}^c_S(\Gamma_R)}(P)\rightarrow \Aut_\mathcal{L}(P)\right)$. Then $\Theta$ is a signalizer functor for $\Gamma _R$. 
\end{Proposition}
\underline{Proof: } We have to show that
$\Theta(P)$ is a complement of $Z(P)$ in $C_{\Gamma_R}(P)$, and for $g\in\Gamma_R$ and $\mathcal{F}$-centric
subgroups $P,Q\leq S$ with $gPg^{-1}\leq Q$, then $\Theta(Q)\leq
g\Theta(P)g^{-1}$.
 We use notation as if 
$\Gamma_R=\Aut_{\mathcal{L}}(S)\ast\dots\ast_{N_S(R_n)}\Aut_{\mathcal{L}}(R_n)\ast
\dots \ast_{N_S(P_i)} N_S(P_i)\ast\dots$. This
way $\delta_P(p)\in\Gamma_R$ makes sense for every
$\mathcal{F}$-centric $P$ (and is equal to $\delta_S(p)$).
We begin by analyzing how we can write $a\delta_S(p)a^{-1}$, where $p\in P$, $a\delta_R(P)a^{-1}=\delta_R(Q)$ and $P,Q\leq N_S(R)$, $a\in\Aut_\mathcal{L}(R)$. By these assumptions there must be a $q\in Q$ with $\delta_S(q)=\delta_{R}(q)=a\delta_R(p)a^{-1}=a\delta_S(p)a^{-1}$. 
Using results from the proposition above, we see that $\delta_Q(q)=\widehat{\delta_R(q)}=\widehat{a\delta_R(p)a^{-1}}=
\widehat{a}\widehat{\delta_R(p)}\widehat{a^{-1}}=\widehat{a}\delta_P(p)\widehat{a}^{-1}$.
From the third axiom of linking systems it follows that
$\delta_Q(\pi(\widehat{a})(p))=\widehat{a}\delta_P(p)\widehat{a}^{-1}=\delta_Q(q)$.
Thus in $\Gamma_R$ we have
 $a\delta_S(p)a^{-1}=a\delta_R(p)a^{-1}=\delta_R(q)=\delta_{Q}(q)=\delta_{Q}(\pi(\widehat{a})(p))=\delta_S(\pi(\widehat{a})(p))$.
Using this we can show that $\Theta(P)\leq C_{\Gamma_R}(P)$. So let
$g\in\Theta(P)$, ie. $g\in N_{\Gamma_R}(P)$ with $\psi_P(g)=\id$. We
have a decomposition of $g$ as $g=a_1\dots a_n$ and $P_i$ etc. as
discussed at the beginning of this section. For every $p\in P$ we have:
$g\delta_S(p)g^{-1}=a_1\dots a_n\delta_S(p) a_n^{-1}\dots a_1^{-1}=a_1\dots a_{n-1}\delta_S(\pi(\widehat{a_n})(p)) a_{n-1}^{-1}\dots a_1^{-1}=\dots=\delta_S(\pi(\widehat{a_1})\circ\dots\circ\pi(\widehat{a_n})(p))=\delta_S(\pi(\widehat{a_1}\circ\dots\circ \widehat{a_n})(p))=\delta_S(\pi(\id)(p))=\delta_S(p)$.
We need to show that $\Theta(P)$ is a complement of $Z(P)$ in
$C_{\Gamma_R}(P)$. For this it suffices to show that $\psi_P|_{Z(P)}$ is an isomorphism onto the image of $\psi_P|_{C_{\Gamma_R}(P)}$.
If $z\in Z(P)$, then $\widehat{\delta_S(z)}=\delta_P(z)$, so $\psi_P$ is injective on $Z(P)$ because $\delta$ is.
For surjectivity, assume that $g=a_1\dots a_n\in C_{\Gamma_R}(P)$.
Then for every $p\in P$ we have $\delta_S(p)=g\delta_S(p)g^{-1}$, so
by our prior work we see that
$\delta_S(p)=\delta_S(\pi(\widehat{a_1}\circ\dots\circ\widehat{a_n})(p))$,
so we have $\id=\pi(\widehat{a_1}\circ\dots\circ\widehat{a_n})$. The
Kernel of $\pi$ is $\delta_P(Z(P))$, so there is a $z\in Z(P)$ with
$\delta_P(z)=\widehat{a_1}\circ\dots\circ\widehat{a_n}$. This means
that
$\psi_P(g)=\widehat{a_1}\circ\dots\circ\widehat{a_n}=\delta_P(z)=\widehat{\delta_S(z)}=\psi_P(\delta_S(z))$.
This finishes the proof of the first property.
Assume that $g=a_1\dots a_n\in \Gamma_R$, $P,Q$ are $\mathcal{F}$-centric and $gPg^{-1}\leq Q$. We need to show that $\Theta(Q)\leq g\Theta(P)g^{-1}$. So let $b=b_1\dots b_m \in \Theta(Q)$. We are to show that $g^{-1}bg\in \Theta(P)$. 
Note that we can restrict the $\widehat{b_i}$ used to obtain $\psi_Q(b)$ in order to get composable morphisms yielding $\psi_{gPg^{-1}}(b)$. It follows that as $\psi_Q(b)=\id$, also $\psi_{gPg^{-1}}(b)=\id$.
Thus $\psi_P(g^{-1}bg)=\widehat{a_n}^{-1}\circ\dots\circ \widehat{a_1}^{-1}\circ\widehat{b_1}\circ\dots\circ\widehat{b_m}\circ\widehat{a_1}\circ\dots\circ\widehat{a_n}=\widehat{a_n}^{-1}\circ\dots\circ \widehat{a_1}^{-1}\circ\psi_{gPg^{-1}}(b)\circ\widehat{a_1}\circ\dots\circ\widehat{a_n}=\widehat{a_n}^{-1}\circ\dots\circ \widehat{a_1}^{-1}\circ\widehat{a_1}\circ\dots\circ\widehat{a_n}=\id$.
This concludes the proof that $\Theta$ is a signalizer functor. $\Box$\\[0.3cm]
In the special case that our
$p$-local finite group is the associated one to a finite group $G$ and under the extra assumption that the signalizer functor for $G$ is
trivial, ie. for all centric subgroups $P$ we have $C_G(P)\leq P$ 
the Robinson model has the following shape:
\[\Gamma_R=N_G(S)\ast_{N_S(R_1)}N_G(R_1)\ast\dots\ast_{N_S(R_n)}N_G(R_n).\]

Let $a\in \Aut_\mathcal{L}(R)\cong N_G(R)$ and $P<S$ be $\mathcal{F}$-centric. Assume that $aPa^{-1}=Q<S$. Then $a\in N_G(\langle R,P\rangle,\langle R,Q\rangle)=\Mor_\mathcal{L}(\langle R,P\rangle,\langle R,Q\rangle)$. We can see that the way we identify morphisms in $\mathcal{L}$ with elements in $G$ identifies $\widetilde{a}$ with $a$ again. This happens also after restricting, so $\widehat{a}$ is actually again $a$, only that $a$ was considered as an element in $N_G(R)$ and $\widehat{a}$ is considered an element in $N_G(P,Q)$. Thus, all $\psi_P$ are restrictions of the map $\Gamma_R\rightarrow G$ given by the inclusions, to $N_{\Gamma_R}(P)$. This means that $\Theta(P)=N_{\Gamma_R}(P)\cap \Ker(\Gamma_R\rightarrow G)$.
\subsection{Signalizer functor for $\pi _{K,T}$}
The previous argument can be adapted to the group model $\pi _{K,T}$. 
The crucial observation is that $\delta _{P}$ labels elements in $Aut_{\mathcal{L}}(P)$ and conjugation via $a_i$ never moves elements out of $S$. Then $\Theta$ is a signalizer functor for $\pi _{K,T}$. The signalizer functor for $\pi _K$ can be also computed from the signalizer functor for $\Gamma _K$ since $\pi _K$ is a quotient group of $\Gamma _K$. We omit the details to the interested reader. 

\section{Euler-Poincar\'e characteristic}
We give formulas for the Euler-Poincar\'e Characteristic for the classifying spaces of groups that realize a prescribed fusion system.
Before we state our results we recall basic properties about the Euler characteristic of finite groups and amalgamated products.
\begin{enumerate}
\item If $A$ is free of rank $r$ (the possiblity $r=0$ is allowed), then $\chi (A)=1-r$.
\item If $A$ has a free subgroup $B$ of index $s$, then $\chi (A)=\frac{\chi (B)}{s}$.
\item If $A=B\underset{C}{*}D$ for groups $B$, $D$ with a common subgroup $C$, then $\chi (A)=\chi (B)+\chi (D)-\chi (C)$ whenever the three rightmost quantities are defined.
\end{enumerate}
Notice that $\chi (A)>0$ if and only if $A$ is finite, and $\chi (A)=0$ if and only if $A$ has an infinite cyclic subgroup of finite index.

\begin{Theorem}
Let $\mathcal{F}$ be a fusion system over the finite $p-$group $S$
and $\pi_{LS}$ a model of Leary-Stancu type for $\mathcal{F}$. Then we have
the following formula
\begin{eqnarray*}
\chi (\pi_{LS})=\frac{1}{|S|}-(\underset{1\leq i\leq n}{\sum
}{\frac{1}{|P_i|}})
\end{eqnarray*}
for the Euler characteristic of $\pi_{LS}$.
\end{Theorem}
\underline{Proof:} A model for the classifying space $B\pi_{LS}$ can
be obtained as gluing the spaces $BP_i \times I$ to a model for the
space $BS$. The Mayer-Vietoris-Sequence for the classifying space of
this group proves the statement.$\Box$
\begin{Theorem} Let $\mathcal{F}$ be a saturated fusion system over the finite $p-$group $S$ and $\pi _{K,T}$ and $\pi_{LS}'$ be group models for $\mathcal{F}$ as above. Then
  \[\chi(\pi_{K,T})=\sum\limits_{i=1}^{n}\frac{1}{|K_i|}-\sum\limits_{(i,j)\text{ with }K_{i,j}\in T}\frac{1}{|K_{i,j}|}\text{ and }
\chi(\pi_{LS}')=\frac{1}{|S|}+\sum\limits_{i=1}^n \frac{1-\ord(\phi_i)}{|P_i|\cdot |\ord(\phi_i)|}.\]

for the Euler characteristic of $\pi_{K,T}$. Moreover we have
\begin{eqnarray*}
\chi (\pi_{K,T})=\frac{d_{\mathcal{F}}}{|S|lcm\{K_i:K_{1,i},1 \leq i\leq
n\}}
\end{eqnarray*}
for some negative integer $d_{\mathcal{F}}$.
\end{Theorem}
\underline{Proof:}
The statement about the Euler characteristic is
just a fact about the Euler characteristic of amalgams and of finite
groups \cite{Serre}. The Euler characteristic of $\pi_{LS}'$ satisfies
\[\chi(\pi_{LS}')=\frac{1}{|S|}+\sum\limits_{i=1}^n \frac{1}{|\langle P_i,t_i \rangle|}-\frac{1}{|P_i|}=\frac{1}{|S|}+\sum\limits_{i=1}^n \frac{1}{|P_i|\cdot|\ord(\phi_i)|}-\frac{1}{|P_i|}=\frac{1}{|S|}+\sum\limits_{i=1}^n \frac{1-\ord(\phi_i)}{|P_i|\cdot |\ord(\phi_i)|}.\Box \]
\begin{Theorem}
Let  $G=S\wr
\Sigma _{e(X)}$ constructed by Park \cite{Sejong}. Then we have $\chi (G)=\frac{1}{|S|\dot |e(X)|\dot |e(X)|!}$.
\end{Theorem}
\underline{Proof:} For a finite group $G$ we have  $\chi (G)=\frac{1}{|G|}$. So in this case we obtain $\chi (G)=\frac{1}{|G|}=\frac{1}{|(S\times\cdots\times S )\Sigma_{e(X)}|}=\frac{1}{|S|\dot |e(X)|\dot |e(X)|!}\Box$.
\section{Examples}
We illustrate our results with examples for group models, explicit signalizer functors, applications to the fundamental group of $|\mathcal{L}|$ and $|\mathcal{L}|\pcom$, existence of centric linking systems via group models.
\subsection{Comparison of $\Gamma_R$ and $\pi_R$}
We compare the two group models $\Gamma_R$ and $\pi_R$ taken over the collection of
$\mathcal{F}$-centric-radical or $\mathcal{F}$-essential subgroups. For the
fusion system of $ASL(2,3)$ at the prime $3$ the two models differ,
with $\Gamma_R$ an infinite group, while $\pi_R = ASL(2,3)$. 
%
\begin{Proposition}
\label{pi_eq_gamma}
 Let $\mathcal{F}$ be a saturated fusion system over the finite $p$-group $S$ with
 $\Out_\mathcal{F}(S)=1$ or $S$ abelian. Then $\Gamma_R=\pi_R$.
\end{Proposition}
\underline{Proof:}
 Let $R_1,\dots R_n$ be the $\mathcal{F}-$centric-radical proper subgroups of $S$. 
We show $
  \Aut_\mathcal{L}(R_i<S)\cong N_S(R_i)$. By
  definition $\Aut_\mathcal{L}(R_i<S)$ is the preimage in $\Aut_\mathcal{L}(S)$ of the subgroup
   $\Aut_\mathcal{F}(R_i<S)$ of $\Aut_\mathcal{F}(S)$, which leaves $R_i$ invariant. So in this case we
   have $\Aut_\mathcal{L}(S)=L_S\cong S$ and $\Aut_\mathcal{F}(S)\cong S/Z(S)$. 
   The
   map $\Aut_\mathcal{L}(S)\rightarrow \Aut_\mathcal{F}(S)$ is the canonical one $S\rightarrow S/Z(S)$.
   Conjugation by an element of $S$ leaves $R_i$ invariant iff that element is in $N_S(R_i)$. So
   $\Aut_\mathcal{F}(R_i<S)\cong N_S(R_i)/Z(S)$. 
   Taking the
    preimage, we see  $\Aut_\mathcal{L}(R_i<S)$ is generated by $Z(S)$ and $N_S(R_i)$, so as
    $Z(S)<N_S(R_i)$ we have $\Aut_\mathcal{L}(R_i<S)\cong N_S(R_i)$. Thus, 
    $\pi_R$ and $\Gamma_R$ coincide. If $S$ is abelian for no proper subgroup $R<S$ we have $C_S(R)=S$ is not contained in
 $R$ there are no proper $\mathcal{F}-$centric-radical subgroups. So $\pi_R=L_S=\Gamma_R$. $\Box$
%
%
\begin{Proposition}
\label{pi_neq_gamma}
Let $(S,\mathcal{F},\mathcal{L})$ be the $3$-local finite group associated with $G=ASL(2,3)$. Then $\pi_R\neq \Gamma_R$, where the models are taken over the $\mathcal{F}$-essential subgroups.
\end{Proposition}
\underline{Proof:}
Recall $G=C_3^2\rtimes SL(2,3)$ and $S$ is an extraspecial group of order $3^3$ and exponent $3$
 generated by the elements $((1,0)^t,I)$, $((0,1)^t,I)$, where $I$ is the identity matrix and
$((0,0),D)$ with $D=\left(\begin{array}{cc} -1 & -1 \\
1&0\end{array}\right)$.There is a unique proper $\mathcal{F}-$centric-radical subgroup denoted $R$. We claim that $R$ is the normal
subgroup $C_3^2$ of $G$. We can first see that $C_G(R)=R$ so
$C_G(R)=Z(R)$, and in particular $R$ is $3$-centric in $G$. As $R$ is normal in $S$, it is fully normalized. 
So $R$
 is $\mathcal{F}$-centric. 
As $R$ is abelian, $\Inn(R)=1$. So
 $\Out_\mathcal{F}(R)=\Aut_\mathcal{F}(R)=N_G(R)/C_G(R)=G/R\cong
 SL_2(3)$. This group has a strongly $3$-embedded subgroup, so $R$ is $\mathcal{F}$-essential and in particular $\mathcal{F}$-radical.
Since $C_G(R)<S$ we have $C_G(S)<S$. So $C_G(S)=Z(S)$ and thus $L_S=T_G(S,S)=N_G(S)=<S,\left((0,0)^t,-I\right)>$.
For $\Gamma_R$ we amalgamate over $N_S(R)=S$. So $\Gamma_R=L_S\ast_S L_R=<S,\left((0,0),-I\right)>\ast_S G$.
Denote the copy of $\left((0,0)^t,-I\right)$ in $G$ by $t$ and the copy in the $L_S$ by $t'$. Then $t't$ is a cyclically reduced element, and thus is of infinite order \cite[page 5]{Serre}.
However, for $\pi_R$ we amalgamate over $\Aut_\mathcal{L}(R<S)$. 
The restriction of conjugation by an element of $N_G(S)$ to $R$ is also an automorphism, as $N_G(R)=G$. So $\Aut_\mathcal{F}(R<S)=\Aut_\mathcal{F}(S)$ and the preimage is the full group, $\Aut_\mathcal{L}(R<S)=\Aut_\mathcal{L}(S)=L_S$. So $\pi_R=L_S\ast_{L_S} L_R\cong L_R=G$. Thus in this case $\pi_R$ is finite while $\Gamma_R$ is infinite. $\Box$
\subsection{Examples of signalizer functors and the fundamental group}
We describe  signalizer functors for various group models by giving generators for the groups $\Theta(P)$ and illustrate the relations with the fundamental group. It is of particular interest in relation to the eighth problem of Oliver's list \cite{AKO} that we can relate group models for fusion systems and the fundamental group since there is no structural pattern for the fundamental group known so far.
\begin{Proposition}  For all fusion systems over $D_8$ we have $\Gamma _R=\pi _R $. For an appropriate choice of $\mathcal{F}^c$-centric subgroups $\Gamma _R=\pi _R \cong \pi _1(|\mathcal{L}|)$, and the signalizer functor for these models is trivial. For the trivial fusion system over $D_8$ we have $\Gamma _R=\pi _R =D_8\cong \pi _1(|\mathcal{L}|)$ and $\pi _1(|\mathcal{L}^c_{D_8}(D_8)|^\wedge_2 )\cong D
_8$. For $\Gamma _R=\pi _R =\Sigma _4\cong \pi _1(|\mathcal{L}|)$ a group model for the $2$-local finite group associated with $\Sigma _4$ we obtain $\pi _1(|\mathcal{L}^c_{D_8}(\Sigma _4)|^\wedge_2 )\cong C_2$.
For $\Gamma _R=\pi _R =\Sigma _4\underset{D_8}{*}\Sigma _4\cong \pi _1(|\mathcal{L}|)$ a model for the $2$-local finite group associated with $PSL_2(7)$ we have $\pi _1(|\mathcal{L}^c_{D_8}(PSL_2(7))|^\wedge_2 )\cong 1$.

\end{Proposition}
\underline{Proof:} It follows from Proposition \ref{pi_eq_gamma} that for all fusion systems over $D_8$ we have $\Gamma _R=\pi _R$. The signalizer functor for the trivial fusion system over $D_8$ is trivial since $D_8$ is a $2-$group and $\pi _1(|\mathcal{L}^c_{D_8}(D_8)|^\wedge_2 )\cong (BD_8)^\wedge_2\cong D
_8$. For the 
$2$-fusion system of $\Sigma _4$ we have $\pi _1(|\mathcal{L}^c_{D_8}(\Sigma _4)|^\wedge_2 )\cong (\Sigma _4)^\wedge_2\cong C_2$. Recall that $ PSL_2(7)$ has a 2-Sylow subgroup isomorphic to the dihedral group of eight elements $D_8$ and denote by 
$V$ and $W$ two representatives of the respective conjugacy classes of the $\mathcal{F}$-centric-radical subgroups of $D _8$. We know \cite[Example 8.8]{IntroMarkus} that $ Aut_{\mathcal{F}}( V ) = Aut_{\mathcal{F}} ( W ) = \Sigma _3$ and therefore we have that $Aut_{\mathcal{L}}( V ) = Aut_{\mathcal{L}}( W ) = \Sigma _4$ and
 $\mathcal{G}=\Sigma _4\underset{D_8}{*}\Sigma _4 $ is a model of Robinson type for the 2-local finite group associated to $PSL(2,7)$. Since $N_{\mathcal{G}}(V)=N_{\mathcal{G}}(W)=\Sigma _4$ and $N_{\mathcal{G}}(D_8)=D_8$ we have $\Theta (V)=\Theta (W)=\Theta (S)=1$. It follows from [3, Theorem 3.5] that
$B\mathcal{G}$ is weakly equivalent to $|\mathcal{L}^c_{D_8}(PSL_2(7))|$  and therefore we obtain $\pi _1(|\mathcal{L}^c_{D_8}(PSL_2(7))|^\wedge_2 )\cong \pi _1(B\mathcal{G}^\wedge_2)\cong  (\Sigma _4\underset{D_8}{*}\Sigma _4)^\wedge_2\cong1$. $\Box $
\begin{Proposition}
 Let $(S,\mathcal{F},\mathcal{L})$ be the $3$-local finite group associated with $G=ASL(2,3)$ and $S,R$ as in Proposition \ref{pi_neq_gamma}. The four proper $\mathcal{F}$-centric subgroups of $S$ are $R$ and $C_{00},C_{10},C_{01}$, where $C_{ij}=\langle((1,1)^t,I),((i,j)^t,D)\rangle$, and form three $\mathcal{F}$-conjugacy classes.  The groups $C_{10}$ and $C_{01}$ are $\mathcal{F}$-conjugate. Furthermore $\Gamma_R=C_3^2\rtimes\langle D,-\id\rangle \ast_S G$.  A signalizer functor $\Theta$ for $\Gamma_R$ is given by $\Theta(S)=\Theta(C_{00})=\Theta(C_{10})=\Theta(C_{01})=\langle tt' \rangle$. The group $\Theta(R)$ will be described below. The signalizer functor for $\pi _R$ is trivial.
\end{Proposition}
\underline{Proof:} 
There are 4 rank 2 elementary abelian $\mathcal{F}$-centric proper subgroups of $S$  \cite[Lemma 3.2]{RV}. Since $D(1,1)^t=(1,1)^t$, we have that $Z(S)=\langle\left((1,1)^t,I\right)\rangle$. Any element of the set $R-Z(S)$ and $Z(S)$ generate $R$. We are reduced to subgroups generated by $Z(S)$ and an element of the form $((i,j)^t,D)$. Because $\left((1,1)^t,I\right)\in Z(S)$ we need only consider $(i,j)\in\{(0,0),(1,0),(0,1)\}$. This proves the first part of the statement.
The group $S$ is its own $\mathcal{F}$-conjugacy class. 
As $R\vartriangleleft G$, $R$ also is its own $\mathcal{F}$-conjugacy class. It remains to consider the $\mathcal{F}$-centric subgroups of the form $C_{ij}$. If two of them are $\mathcal{F}$-conjugate via $((a_1,a_2)^t,A)\in G$, then $A\in N_H(\langle D\rangle)=C_H(\langle D\rangle)$ with $H=SL(2,3)$. Direct calculation shows
$((a_1,a_2)^t,A)((b_1,b_2)^t,D)(-A^{-1}(a_1,a_2)^t,A^{-1})
=((a_2-a_1,a_2-a_1)^t,I)\cdot (A(b_1,b_2)^t,D)$.
%
We can read off that $C_{00}$ is $\mathcal{F}$-conjugated only to itself. As $(-I(1,0)^t,D)=((-1,-1)^t,I)\cdot ((0,1),D)$ we also obtain that $C_{10}$ and $C_{01}$ are $\mathcal{F}$-conjugate via $((0,0)^t,-I)$.
Recall the Robinson model for $G$ is
$\Gamma_R=\langle S,t'\rangle\ast_S G$
where $t'=((0,0)^t,-I)$. Denote by $t$ the respective element in the $G$ part of the amalgam.
 As $S$ is its own conjugacy class, we have $N_\Gamma(S)=N_{L_S}(S)\ast_S N_{L_R}(S)=\langle S, t'\rangle \ast_S \langle S,t\rangle$.
All elements in this amalgam can be written as a reduced word \cite{Serre}, which in this case where $\langle S, t\rangle$ has only two classes mod $S$ leads to elements of the form $st'^\alpha tt't\dots t'^\beta$ where $\alpha,\beta\in\{0,1\}$ and $s\in S$. Multiplying this as an element in $G$ we obtain  $s\cdot(0,-I)^m$, where $m$ is the number of characters after the $s$ in the word  $st'^\alpha tt't\dots t'^\beta$. This is $(0,I)$ iff $s=(0,I)$ and $m=2n,n\geq 0$. Thus as $(tt')^{-1}=t't$, we conclude $\Theta(S)=\langle tt' \rangle$.
The group $H\in Syl_3(\langle D\rangle)$ has a normal
complement generated by
 $
E_1=\left(\begin{array}{ccc}
1 & 1  \\
1 & -1 \\
          \end{array}
\right) \text{, }
E_2=\left(\begin{array}{ccc}
-1 & 1  \\
1 & 1 \\
          \end{array}
\right)$, $-I$.
We have $N_\Gamma(R)=N_{L_S}(R)\ast_SN_{L_R}(R)=L_S\ast_SL_R=\Gamma$.
Again, reduced words have the form
$((s_1,s_2)^t,D^\lambda)t'^\alpha(0,A_1)t'(0,A_2)t'\dots t'(0,A_n)t'^\beta$.
Here $\alpha,\beta\in\{0,1\}$, $\lambda\in\{0,1,2\}$ and $A_i\in\langle -I,E_1,E_2\rangle$.
Here we consider all the elements except the $t'$ as elements in $L_R$. If we multiply this
together in $G$ we get $((s_1,s_2)^t,D^\lambda\cdot (-\id)^{\alpha+\beta+n-1}A_1\dots A_n)$.
If this is to be $(0,I)$, then certainly $s_1=s_2=0$, $\lambda=0$. It remains to determine for which $A_1,\dots,A_n$ the product is $I$ or $-I$. As $tt'\in \Theta(R)$ is already clear, we can reduce to the case $A_1\dots A_n=I$,
where $A_i\in\langle -I,E_1,E_2\rangle$. 
Consider this latter group as a monoid, and consider the presentation with generators $-I, E_1,E_2$. Let $V$ be a set of generators of the relations. An element of $V$ is of the form $B_1\dots B_n$ with $B_i\in\{-I,E_1,E_2\}$, and clearly all the elements $t'^{\alpha_0} B_1t'^{\alpha_1}\dots B_nt'^{\alpha_n}$ are in $\Theta(R)$, where $\alpha_i\in\{0,1\}$ such that $\Sigma \alpha_i$ is even. Conversely, if we have $A_1,\dots,A_n$ with $A_i\in\langle -I,E_1,E_2\rangle$ and such that the product is $I$, then first write $A_i=B_{i1}\dots B_{il_i}$ with $B_i\in\{-I,E_1,E_2\}$. Then $t'^\alpha A_1 t'\dots t' A_nt'^\beta=t'^\alpha B_{11}\dots B_{1l_1}t'\dots t' B_{n1}\dots B_{nl_n}t'^\beta=t'^\alpha B_{11}\dots t'B_{1l_1}t'\dots t' B_{n1}\dots B_{nl_n}t'^\beta$.
The product $B_{11}\dots B_{nl_n}$ is $I$, so it can be written as the product of elements in $V$. Thus (recall that $t'^2=I$) we can also write $t'^\alpha A_1 t'\dots t' A_nt'^\beta$ as the product of elements of the form  $t'^{\alpha_0} B_1t'^{\alpha_1}\dots B_nt'^{\alpha_n}$ just discussed. Thus $\Theta(R)$ is generated by these elements and $tt'$. 
We have $N_H(\langle D\rangle)=\langle D, -I\rangle$,
and 
therefore $N_G(C_{0,0})=C_3^2\rtimes \langle D,-I\rangle$ and
$N_\Gamma(C_{00})=N_{L_S}(C_{00})\ast_S N_{L_R}(C_{00})=L_S\ast_S L_S$.
We see that the normalizer is identical to the one of $S$.
Thus, $\Theta(C_{00})=\Theta(S)=\langle tt'\rangle$.
Next, $\langle N_G(C_{10}), N_G(C_{1,0},C_{0,1}), N_G(C_{0,1},C_{1,0}), N_G(,C_{0,1})\rangle = \langle S, t\rangle$. 
Conjugation with $t$ or $t'$ switches $C_{0,1}$ and $C_{1,0}$. 
Thus as above  $N_\Gamma(C_{0,1})=N_\Gamma(C_{1,0})$ contain exactly the elements of the form $stt'\dots tt'$ and $st't\dots t't$.
We have $\Theta(C_{10})=\Theta(C_{01})=\langle tt'\rangle$. \\
The signalizer functor for $\pi _R$ is trivial because $\pi _R=ASL(2,3)$ as proved in Proposition \ref{pi_neq_gamma} and in $ASL(2,3)$ all $\mathcal{F}$-centric subgroups contain their centralizers. $\Box$
\begin{Proposition}
Let $\mathcal{F}$ be a nontrivial saturated fusion system over the finite $p$-group $S$ and $\mathcal{G}$ a simple (possibly infinite) group model for $\mathcal{F}$ with $B\mathcal{G}\simeq |\mathcal{L}|$. Then $\pi _1(|\mathcal{L}|\pcom)$ is trivial.
\end{Proposition}
\underline{Proof:} This follows from Proposition \ref{BG SIMEQ |L|} and the fact that $\mathcal{G}$ is simple since in this case $\pi _1(|\mathcal{L}|\pcom\cong\mathcal{G}/M$ where 
$M$ is the nontrivial maximal $p$-perfect subgroup of $\mathcal{G}$. $\Box$\\
We can show that the fundamental group of the Solomon groups ist trivial because the group model constructed by Libman and the second author \cite{LS} fulfills the conditions of the previous Proposition and the Solomon fusion systems are simple.\\
If $G$ is a finite group then one can detect the $p$-quotients of $G$ from the fusion system of $G$. A $p$-quotient yields a normal subgroup of index $p$, which corresponds to a normal subsystem. If $G$ has a simple fusion system it has no normal subsystems, so has no $p$-quotients. 
\begin{Proposition}
Let $\mathcal{F}$ be a saturated fusion system over the finite $p$-group $S$ and $\mathcal{G}$ be the Leary-Stancu model with the set of morphisms a set of morphisms in the associated centric linking system which fulfill Alperin's theorem. Then we have a surjective map $\mathcal{G}\rightarrow \pi _1(|\mathcal{L}|)$.
\end{Proposition}
\underline{Proof:} This follows from Alperin's fusion theorem for linking systems and the construction of the map $r$ in this case. $\Box$\\[0.3cm]
Since $B\mathcal{G}$ is $p-$bad in general as we will show below this does not allow to compute the fundamental group of $|\mathcal{L}|\pcom$.
\subsection{Existence of classifying spaces via group models for fusion systems}
The group model $\pi_R$ taken over the collection of $\mathcal{F}$-centric-radical subgroups  provides a solution to the existence of centric linking systems for saturated fusion systems over $p$-groups of order $p^3$. 
\begin{Proposition}
Let $\mathcal{F}$ be a saturated fusion system over the finite $p$-group $S$ and $|S|=p^3$. Then $B\pi _R\simeq |\mathcal{L}|$.
\end{Proposition}
\underline{Proof:} Let $R_1,\cdots, R_n$ be representatives of proper $\mathcal{F}$-centric-radical subgroups $S$. Since $|R_i|=p^2$ we have $N_S(R_i)=S$ for all $i=1, \cdots ,n$. By \cite[Theorem 3.5]{BCGLO1} we have $B\pi _R\simeq |\mathcal{L}|$.$\Box$\\

Note that this in particular applies to the exotic Ruiz-Viruel examples \cite{RV} and therefore allows an easy computations of the fundamental groups which is trivial in all three cases since they are all simple fusion systems.\\

The following Proposition illustrates that the homotopy type
 of the group model $\pi _{K,T}$ changes depending on the which groups one amalgamates over in contrast with \cite[Theorem 3.5]{BCGLO1}. It is however not a proof that $\pi_{K,T}$ does not provide a solution to the existence conjecture. This remains open.
\begin{Proposition}
 Consider the $2$-local finite group of $C_2^3\rtimes \GL(3,2)$.
 Then $H^*(\pi_R,\F_2)\ncong H^*(\mathcal{F})$.
\end{Proposition}
\underline{Proof:}
We prove that $H^2(\pi_R;\F_2)\ncong H^2(C_2^3\rtimes \GL(3,2);\F_2)$. The proof will begin by stating some general facts about the groups involved, and then go on to first determine the $\mathcal{F}$-centric-radical subgroups. Afterwards, individual cohomology groups will be calculated and used to derive a contradiction.
Elements of $G$ are of the form $((a_1,a_2,a_3)^t,A)$ with
$((a_1,a_2,a_3)^t,A)\cdot ((b_1,b_2,b_3)^t,B)=((a_1,a_2,a_3)+A(b_1,b_2,b_3)^t,AB)$.
The group $\GL(3,2)$ 
has 
generators 
$A_2=\left(\begin{array}{ccc}
1 & 1 & 0 \\
0 & 1 & 0 \\
0 & 0 & 1
          \end{array}
\right), A_3=\left(\begin{array}{ccc}
0 & 0 & 1 \\
1 & 0 & 0 \\
0 & 1 & 0
          \end{array}
\right)$
with $A_2^2=1, A_3^3=1, (A_2A_3)^7=1, [A_2,A_3]^4=1$. %
A Sylow $2$-subgroup $T$ of $\GL(3,2)$ is given by the subgroup of
upper triangular matrices in $\GL(3,2)$ 
with generators $A_2$ and
$A_4=\left(\begin{array}{ccc}
1 & 1 & 1 \\
0 & 1 & 1 \\
0 & 0 & 1
          \end{array}
\right)$
satisfying $A_2^2=1, A_4^4=1, A_2A_4A_2A_4=1$.We have $C_2^3\rtimes T=S\in Syl_2(G)$.
\begin{Lemma}
 Every centric-radical subgroup of $S$ contains $C_2^3$, no subgroup of $S$ of order $16$ is $\mathcal{F}$-radical, and the group $C_2^3\rtimes \langle A_4\rangle$ is not
$\mathcal{F}$-radical.
\end{Lemma}
\underline{Proof:}
A simple calculation shows that $\left((1,0,0)^t,\id\right)\in Z(S)$.
Assume $R$ to be centric-radical. Thus $N_G(R)/C_G(R)\cdot R$ has no
proper normal $2$-subgroup. Let $Z:=C_2^3\cap N_G(R)$. Then as
$C_2^3 \vartriangleleft G$, we see that $Z/C_G(R)\cdot R$ is a
normal $2$-subgroup of $N_G(R)/C_G(R)\cdot R$. 
Assume there is a $z=((z_1, z_2, z_3),\id)\notin R$, and chose this
element so that $z_3=1$ only if all elements of this form with
$z_3=0$ are already in $R$.
We wish to show that this element is nontrivial
in $Z/C_G(R)\cdot R$, contradicting the fact that $R$ is radical.
 As $R$ is centric, $C_S(R)<R$.  As
$z\notin R$, but $z\in S$, we obtain $z\notin C_G(R)\cdot R$. It
remains to show that $z\in N_G(R)$. This can be done by calculation:
\[((z_1, z_2, z_3)^t,\id)\cdot\left((a_1, a_2, a_3)^t,\left(\begin{array}{ccc}
1 & a_{12} & a_{13} \\
0 & 1 & a_{23} \\
0 & 0 & 1
          \end{array}
\right)\right) \cdot((z_1, z_2, z_3)^t,\id)\]
 \[=((a_{12}z_2+a_{13}z_3,a_{23}z_3,0),\id)\cdot\left((a_1,a_2,a_3),\left(\begin{array}{ccc}
 1 & a_{12} & a_{13} \\
 0 & 1 & a_{23} \\
 0 & 0 & 1
           \end{array}
\right)\right).\]
By our assumptions on  the value of $z_3$ we can conclude
that this is again an element in $N_G(R)$.

Thus the groups that remain as potential centric-radical subgroups of $S$
must have oder $8,16,32$ or $64$ and have the form $C_2^3\rtimes H$
with $H<T$. 
Such a subgroup would have the form $C_2^3\rtimes \langle A \rangle$
with $A\in T$ of order 2. Since $\GL(3,2)\cong
\PSL(3,2)\cong \PSL(2,7)=LF(2,7)$ is a CIT-group \cite{Suz61}, as $7$ is Mersenne
prime, which means that the centralizer of any
involution is a $2$-group. In our case we obtain:
$N_G(C_2^3\rtimes \langle A \rangle)=C_2^3\rtimes N_{\GL(3,2)}(\langle A \rangle)=C_2^3\rtimes C_{\GL(3,2)}(A)$
is a $2$-group, too. From this we can conclude that no such subgroup can be $\mathcal{F}$-radical. 
The group $ \langle A_4\rangle$ has exactly one element of order $2$
 so $N_{\GL(3,2)}(\langle
A_4\rangle)<C_{\GL(3,2)}(\langle A_4^2\rangle)$, which is a $2$-group as $\GL(3,2)$
is a CIT-group. This means that $N_G(C_2^3\rtimes \langle
A_4\rangle)$ is a $2$-group, and thus $C_2^3\rtimes \langle
A_4\rangle$ is not $\mathcal{F}$-radical. $\Box$

We compute normalizers and centralizers of
$C_2^3\rtimes H_i$ and prove that they are
$\mathcal{F}$-centric-radical.
\begin{Lemma}
There are four subgroups of $S$ left, which
are not ruled out as centric-radical by the previous arguments: $R_1=C_2^3\rtimes H_1$, $R_2=C_2^3\rtimes H_2$, $R_3=C_2^3\rtimes H_3$, $R_4=C_2^3\rtimes H_4=S$. These are indeed centric-radical.
\end{Lemma} 
\underline{Proof:}
We need
only consider subgroups of the form $C_2^3\rtimes H$ with $H<
T\cong D_8$, and $|H|\neq 2$, $H\neq \langle A_4\rangle$. There are four such subgroups, $H_1=\{\id\},H_2=\langle A_2,A_4^2\rangle,H_3=\langle A_2A_4,A_4^2\rangle,H_4=T$
with $|H_1|=1, |H_2|=4, |H_3|=4, |H_4|=8$. 
Let $A_z:=A_4^2$ and $A_2'=A_2A_4$.
Suppose that $\left((z_1,z_2,z_3),Z\right)\in C_G(R_1)$. 
A calculation shows that $Z(a_1,a_2,a_3)^t=(a_1,a_2,a_3)^t$ for all $(a_1,a_2,a_3)^t\in C_2^3$. So we can conclude that $Z=\id$
must hold. Since $C_2^3$ is abelian, we have
$C_G(R_1)=R_1$, so $R_1$ is $\mathcal{F}$-centric.
$N_G(R_1)/R_1\cdot C_G(R_1)=G/R_1\cong \GL(3,2)$. This has no proper
normal $2$-subgroups. Indeed, one Sylow subgroup is given by the
upper triangular matrices, another by the lower triangular matrices,
and their intersection is trivial.
Thus $R_1$ is $\mathcal{F}$-centric-radical. 
The group $S$ is always centric-radical. We claim $N_G(S)=S$. As
$N_G(S)=N_G(C_2^3\rtimes T)=C_2^3\rtimes N_{\GL(3,2)}(T)$, we only
need that $N_{\GL(3,2)}(T)=T$. There are only two elements in
$T\cong D_8$ of order $4$, and they have the same square, an element
of order $2$, $A_z$. So any element in $\GL(3,2)$ normalizing $T$
must leave $A_z$ fixed, ie. $N_{\GL(3,2)}(T)<C_{\GL(3,2)}(A_z)$,
which is a $2$-group. So we must have $N_{\GL(3,2)}(T)=T$. With
$N_G(S)=S$, we have that $C_G(S)<S$ and thus $C_G'(S)=1$. Thus,
$\Aut_\mathcal{L}(R_4)=S/1=S$. 
 In the following statements containing $i$ are meant for both $i=2$ and $i=3$.
Let
$A_{3r}=\left(\begin{array}{ccc}
1 & 0 & 0 \\
0 & 0 & 1 \\
0 & 1 & 1
          \end{array}
\right),A_{3l}=\left(\begin{array}{ccc}
1 & 1 & 0 \\
1 & 0 & 0 \\
0 & 0 & 1
          \end{array}
\right).$
These are elements of order $3$. We claim that
$N_{\GL(3,2)}(H_2)=\langle T,A_{3r} \rangle$ and
$N_{\GL(3,2)}(H_3)=\langle T,A_{3l} \rangle$.
Firstly, $H_i \vartriangleleft T$, as it is of index $2$. We can
show by calculation that $A_{3r}$ normalizes $H_2$ and $A_{3l}$
normalizes $H_3$, eg. $A_{3r}\cdot A_z \cdot A_{3r}^{-1}=A_2, A_{3l}\cdot A_z \cdot A_{3l}^{-1}=A_z\cdot A_2'$. Thus 
the normalizers contain the subgroups given above. The group already given must have order at
least $24$, as $|T|=8$ and $A_{3l}$ and $A_{3r}$ are of order $3$.
Thus the normalizers are either those claimed or the whole group $\GL(3,\F_2)$.
Direct calculation shows that $A_{3r}$ does not normalize $H_3$ and $A_{3l}$
does not normalize $A_{3l}$, as eg. $A_{3l}\cdot A_2 \cdot$. 
Since $R_1<R_i$, and $R_1$ is $\mathcal{F}$-centric, the $R_i$ are $\mathcal{F}$-centric, too. Thus
$|N_G(R_i)/C_G(R_i)\cdot
R_i|=|N_G(R_i)|/|R_i|=\frac{8\cdot 24}{8\cdot 4}=6$. The group
$N_G(R_i)/C_G(R_i)\cdot R_i$ is thus either $S_3$ or $C_6$. For it
to be $2$-reduced, we need it to be $S_3$. It suffices to
show that it is not abelian.
For $i=2$: $A_{3r}^{-1}\cdot A_4\cdot A_{3r}\cdot A_4^{-1}\notin S$,
for $i=3$: $A_{3l}^{-1}\cdot A_4\cdot A_{3l}\cdot A_4^{-1}\notin S$.
This shows that the above group is $S_3$, which has no normal $2$-subgroup.
Thus $R_2$ and $R_3$ are $\mathcal{F}$-radical.$\Box$\\[0.3cm]
The above calculations show that $C_G'(R_i)=1$ for $i\in\{1,2,3,4\}$, and the following table holds:

\begin{tabular}{c|c|c|c}
 $i$ & $|R_i|$ & $R_i$ & $\Aut_\mathcal{L}(R_i)$ \\\hline
 $1$ & $8$ & $C_2^3\rtimes 1$ & $G=C_2^3\rtimes \GL(3,2)$ \\\hline
 $2$ & $32$ & $C_2^3\rtimes \langle A_z,A_2\rangle$ & $C_2^3\rtimes \langle T, A_{3r}\rangle$ \\\hline
 $3$ & $32$ & $C_2^3\rtimes \langle A_z,A_2'\rangle$ & $C_2^3\rtimes \langle T,A_{3l}\rangle$ \\\hline
 $4$ & $64$ & $C_2^3\rtimes T$ & $S=C_2^3\rtimes T$ \\\hline
\end{tabular}

As $\Aut_\mathcal{L}(S)=S$ is contained in $\Aut_\mathcal{L}(R_i)$ for all $i$, we get:
\[\pi_R=\left(\left(S\ast_S G\right)\ast_S (C_2^3\rtimes \langle T, A_{3r}\rangle)\right)\ast_S (C_2^3\rtimes \langle T,A_{3l}\rangle)\cong G \ast_S \left( C_2^3\rtimes \langle T, A_{3r}\rangle\right)\ast_S (C_2^3\rtimes \langle T,A_{3l}\rangle).\]
%
\begin{Lemma} We have the following results about cohomology groups: $H^1(G)\cong 0$, $H^1(S)\cong C_2^3$, $H^1(C_2^3\rtimes \langle T, A_{3r}\rangle)\cong C_2$, $H^1(C_2^3\rtimes \langle T, A_{3l}\rangle)\cong C_2^2$, and $H^2(G)\leq C_2^2$.
\end{Lemma}
\underline{Proof:}
We have  $H^1(G)\cong\Hom(G,C_2)$. Let $\varphi: G\rightarrow C_2$ be a group homomorphism. All
elements of order coprime to $2$ must be mapped to $0$. The group $\GL(3,2)$
is generated by $A_3$ of order $3$ and $A_2A_3$ of order $7$. Thus,
$(0,A_3)$ and $(0,A_2A_3)$ are elements of order coprime to $2$
which generate $0\rtimes \GL(3,2)$. As there is always a linear
isomorphism sending a specific nonzero element of $C_2^3$ to another
nonzero element, we get that all nonzero elements of the form
$((z_1, z_2, z_3)^t,\id)$ are sent to the same element by $\varphi$,
which thus must be $0$. These elements together with $(0,A_3)$
and $(0,A_2A_3)$ generate $G$, we can conclude that $\varphi=0$, and
so $H^1(G)\cong C_2^0$. 
The derived subgroup of $S=C_2^3\rtimes T$ is the normal closure of
the subgroup generated by commutators of its generators. The group
$S$ is generated by $C_2^3, A_2$ and $A_2'$. 
A calculation shows that $S'$ is the normal closure of the subgroup 
$U$ of elements of the form $((z_1,z_2,0)^t,A_z^\epsilon)$. Since $A_z\leq Z(T)$ we can check by calculation
that $S$ normalizes $U$. %
From all this we know that $S_{ab}=S/S'$ has $64/8=8$ elements. A calculation shows that all nonzero elements have order $2$ in $S_{ab}$. 
It follows that $S_{ab}=C_2^3$, so $H^1(S)\cong
\Hom(C_2^3,C_2)\cong C_2^3$. 
Denotel $L_2=C_2^3\rtimes \langle T, A_{3r}\rangle$. 
Some calculations similar to the case above show that $L_2'$ is the normal closure of
$C_2^3\rtimes \langle A_{3r},A_z\rangle$. The group $\langle A_{3r},A_z\rangle$ has order at least $12$, as it contains a subgroup of order $4$ generated by $A_z$ and $A_{3r}$ of order $3$. It cannot be $\langle T, A_{3r}\rangle$ (of order $24$), as this group is not perfect, and so $\langle A_{3r}.A_z,\rangle$ is a normal subgroup of index $2$.
Thus, $|{L_2}_{ab}|=2$, and so $H^1(L_2)\cong \Hom(C_2,C_2)\cong
C_2$. 
Let $L_3=C_2^3\rtimes \langle T, A_{3l}\rangle$. 
Then $L_3'$ is the normal closure of $(C_2^2\times 0)\rtimes \langle
A_{3l},A_z\rangle$.
Reasoning as above, $\langle A_{3l},A_z\rangle$ has order $12$
and is normal in $\langle T, A_{3l}\rangle$. Since elements in $T$ fix the last coordinate, we conclude that $L_3'=(C_2^2\times
0)\rtimes \langle A_{3l},A_z\rangle$. Thus,
$|{L_3}_{ab}|=\frac{8\cdot 24}{4\cdot 12}=4$. Every element in
${L_3}_{ab}$ is of order $2$, as $A_2^2=\id$.
We conclude 
$H^1(L_3)\cong\Hom({L_3}_{ab},C_2)\cong\Hom(C_2^2,C_2)\cong C_2^2$. 
Since $G=C_2^3\rtimes \GL(3,2)$ using \cite{Tahara} we have $H^2(G,\F_2)\cong H^2(\GL(3,2),\F_2)\oplus \tilde{H}^2(G,\F_2)$
where the following sequence is exact:
$0\rightarrow H^1(\GL(3,2),H^1(C_2^3))\rightarrow \tilde{H}^2(G,\F_2)\rightarrow H^2(C_2^3,\F_2)^{\GL(3,2)}$. By \cite{Griess} and the universal coefficient theorem we have
$H^1(\GL(3,2),\F_2)\cong C_2$.
By direct calculation $H^2(C_2^3,\F_2)^{\GL(3,2)}\cong 0$ so
 we have
$0\rightarrow H^1(\GL(3,2),H^1(C_2^3))\rightarrow \tilde{H}^2(G,\F_2)\rightarrow 0$
and we only need to show $ H^1(\GL(3,2),H^1(C_2^3))\leq C_2$.
Since $H^1(\GL(3,2),H^1(C_2^3))\cong  H^1(\GL(3,2),C_2^3)$
can be
identified with the group we get by dividing the principal crossed
homomorphism from $\GL(3,2)$ to $C_2^3$ out of the crossed
homomorphism.
Let $\varphi$ be a crossed homomorphism. The group $\GL(3,2)$ is generated by
$A_2$ and $A_3$, with the relations given at the beginning. Since
\[0=\varphi(A_2^2)=(A_2+\id)\varphi(A_2)=\left(\begin{array}{ccc}
0 & 1 & 0 \\
0 & 0 & 0 \\
0 & 0 & 0
          \end{array}
\right)\varphi(A_2)\]
 $\varphi(A_2)$ must be of the form $(*,0,*)^t$.
By similar calculations with $A_3^3=\id$ we obtain that the three coordinates of $\varphi(A_3)$ must add up to $0$.
As a crossed homomorphism is already fixed by the image of
generators, we see that there can be a maximum of $4\cdot 4=16$
crossed homomorphisms. There are $8$ different
principal crossed homomorphisms, given by $C_2^3$.
Thus, $|H^2(\GL(3,C_2^3))|\leq 2$. $\Box$\\

\begin{Lemma} $H^2(G)\ncong H^2(\pi_R)$.
\end{Lemma}
\underline{Proof:}
Assume $H^2(G)\cong H^2(\pi_R)$. Then we get a long exact
Mayer-Vietoris sequence:
\[H^1(\pi _R)\rightarrow H^1(G)\times H^1(L_2)\times H^1(L_3)\rightarrow H^1(S)\times H^1(S)\rightarrow H^2(G)\rightarrow \cdots\]
which reduces to $0\rightarrow H^1(\pi _R)\rightarrow C_2^1\times C_2^2\rightarrow
C_2^3\times C_2^3\rightarrow H^2(G)\rightarrow\cdots$. 
The image of the second map has rank at most $3$, so the fourth map must have rank at least $(3+3)-3=3$. However, $H^2(G)$ has rank at most $2$ by the
above calculations which is a contradiction. $\Box$\\[0.3cm]
Recall that the Leary-Stancu model allows to assign a group model to a fusion system in a functorial way. The following proposition however ruins the hope that it could be minimal in almost all cases.
\begin{Proposition}
\label{learystancupbad}
Let $\mathcal{G}$ be a group model for a fusion system $\mathcal{F}$ over the finite $p$-group $S$ associated with a graph of groups $\Gamma$ which is a one-dimensional category and contains at least $2$ loops and such that $B\mathcal{G}\simeq hocolim_{\Gamma}B(-)$. Then the classifying space $B\mathcal{G}$ is $p$-bad for any prime $p$.
\end{Proposition}
\underline{Proof:} It follows from  \cite[Theorem 1.1] {Farjoun}that the space $B\mathcal{G}$ contains a wedge of spheres as a retract which is $p-$bad for all primes \cite{Bousfield}. A space which contains a $p$-bad space as a retract is $p$-bad. $\Box$\\[0.3cm]
In the special case of a Leary-Stancu model we obtain that the classifying space is $p$-bad as soon as the group is more than a single HNN-extension.

Dr.\ Nora Seeliger, Department of Mathematics, University of Haifa, Faculty of Natural Sciences, Science and Education Building, Room 615, Abba Khoushy Avenue 199, 3498838 Haifa, ISRAEL.\\Email:
seeliger@math.univ-paris13.fr


\begin{thebibliography}{breitestes Label}
\bibitem{AschbacherChermak} M. Aschbacher, A. Chermak, \textit{ A group-theoretic approach to a family of 2-local finite groups
constructed by Levi and Oliver}, Ann. of Math. (2), 171 (2010), no.
2, 881--978.
\bibitem{AKO} M. Aschbacher, R. Kessar, B. Oliver, \textit{Fusion systems in algebra and topology}, London Mathematical Society Lecture Note Series: 31, Cambridge University Press, 2011.
\bibitem{Bousfield} A. K. Bousfield, \textit{On the $p-$adic completions of nonnilpotent spaces}, Transactions of the American Mathematical Society, Volume 331, Number 1, May 1992, 335--359.
\bibitem{BK} A. K. Bousfield, D. M. Kan, \textit{Homotopy Limits, Completions and Localizations}, Springer Lecture Notes in Mathematics, Springer-Verlag, Berlin, Heidelberg, New York, 1972.
\bibitem{BCGLO1} C. Broto, N. Castellana, J. Grodal, R. Levi,  B. Oliver, \textit{Subgroup families controlling $p-$local finite groups}, Proc., London Math. Soc. 91 (2005), 325--354.
\bibitem{BLO1} C. Broto, R. Levi B. Oliver, \textit{Homotopy Equivalences of $p-$Completed Classifying Spaces of Finite Groups}, Invent. Math. 151 (2003), 611--664.
\bibitem{BLO2} C. Broto, R. Levi, B. Oliver, \textit{The Homotopy Theory of Fusion Systems}, J. Amer. Math. Soc. \textbf{16} (2003), no. 4, 779--856.
\bibitem{BLO4} C. Broto, R. Levi, B. Oliver, \textit{A Geometric Construction of Saturated Fusion Systems}, Contemp. Math. 399, 2006, 11--39.
\bibitem{Chermak} A. Chermak, \textit{Fusion systems and localities}, preprint 2011.
\bibitem{COS} A. Chermak, B. Oliver, S. Sphectorov, \textit{The linking systems of the Solomon groups are simply connected}, Proc. London Math. Soc., 97 (2208), 209--238.
\bibitem{Farjoun} E. D. Farjoun, \textit{Fundamental group of homotopy colimits}, Advances in Mathematics 182 (2004), 1--27.
\bibitem{Griess} R.\ L.\ Griess, \textit{Schur multipliers of the known finite simple groups}, Bull. Amer. Math. Soc. Volume 78, Number 1 (1972), 68--71.
\bibitem{KanThurston} D. M. Kan, W. P. Thurston, \textit{Every connected space has the homology of a $K(\pi ,1)$}, Topology Vol. 15, pp. 253--258 (1976).
\bibitem{Ian+Radu} I. Leary, R. Stancu, \textit{Realising fusion systems}, Algebra Number Theory 1.1, (2007), 17--34.
\bibitem{LO} R. Levi, B. Oliver, \textit{Construction of 2-local finite groups of a type studied by Solomon and Benson}, Geometry and Topology 6 (2002), 917--990.
\bibitem{LeviOliver} R. Levi, B. Oliver, \textit{Construction of $2$-local finite groups of a type studied by Solomon and Benson}, Geometry and Topology, 6 (2002), 917--990.
\bibitem{LeviOliverC} R. Levi, B. Oliver, \textit{Correction to: Construction of $2$-local finite groups of a type studied by Solomon and Benson}, Geometry and Topology, 9 (2005), 2395--2415.
\bibitem{LS} A. Libman, N. Seeliger, \textit{Homology decompositions and groups inducing fusion systems}, Homology, Homotopy and Applications, Vol. 14 (2012), No. 2, 167--187.

\bibitem{Lib06} A. Libman, \textit{The normaliser decomposition of p-local finite groups}, Algebraic and Geometric Topology, vol 6 (2006), pp. 1267--1288.



\bibitem{IntroMarkus} M. Linckelmann, \textit{Introduction to Fusion Systems}, in: Group Representation Theory, EPFL Press, Lausanne, 2007, pp. 79--113.
\bibitem{MartinoPriddy} J.\ Martino, S.\ Priddy, \textit{Unstable homoropy classification of $BG\pcom$}, Math.\ Proc.\ Camb.\ Phil.\ Soc.\ 137 (2004) 321--347.
\bibitem{Oliver} B.\ Oliver, \textit{Existence and uniqueness of linking systems: Chermak`s proof via obstruction theory}, arXiv:1203.6479.
 \bibitem{MP1} B.\ Oliver, \textit{Equivalences of classifying spaces completed at odd primes}, Math.\ Proc.\ Camb.\ Phil.\ Soc.\ 137, 2004, 321-347.
\bibitem{MP2} B.\ Oliver, \textit{Equivalences of classifying spaces completed at the prime 2}, Mem.\ Amer.\ Math.\ Soc.\ 180 (2006), no.\ 848, vi+102pp.
\bibitem{O4} B.\ Oliver, \textit{ Extensions of linking systems and fusion systems}, Trans. Amer. Math. Soc. 362 (2010), 5483--5500.
\bibitem{Sejong} S.\ Park, \textit{Realizing a fusion system by a single finite group}, Arch. Math. 94 (2010), 405--410.
\bibitem{Quillen} D.\ Quillen, \textit{The spectrum of an equivariant cohomology ring: I}, Annals of Math. 94 (1971), 549--572.
\bibitem{Ragnarsson} K.\ Ragnarsson, \textit{Classifying spectra of saturated fusion systems}, Algebr. Geom. Topol. 6 (2006), 195--252.
\bibitem{Robinson1} G.\ Robinson, \textit{Amalgams, blocks, weights, fusion systems and finite simple groups}, Journal of Algebra 314 (2007), 912--923.
\bibitem{Robinson2} G.\ Robinson, \textit{Reduction (Mod q) of Fusion System Amalgams}, Trans. Amer. Math. Soc. 363 (2011), 1023--1040.
\bibitem{gmffs} N.\ Seeliger, \textit{Group models for fusion systems}, Topology and its Applications 159 (2012), no. 12, 2845--2853.
\bibitem{assclassfisoq} N.\ Seeliger, \textit{Assigning a classifying space to a fusion system up to $F-$isomorphism}, arXiv:1105.3403v2.
\bibitem{RV} A.\ Ruiz, A. Viruel, \textit{Classification of Fusion System over the Extraspecial Group of order $p^3$ and exponent $p$}. Math. Z. 248, 45--65
(2004).
\bibitem{Serre} J.-P. Serre, \textit{Trees}, Springer-Verlag Berlin Heidelberg New York, 1980.

\bibitem{Suz61} M.\ Suzuki, \textit{Finite groups with nilpotent centralizers}, Transactions of the American Mathematical Society 99: 425--470 (1961).
\bibitem{Tahara} K.\ Tahara, \textit{On the Second Cohomology Groups of Semidirect Products}, Mathematische Zeitschrift (1972), Volume 129, Issue 4, pp 365--379.

\end{thebibliography}
\end{document}